\documentclass[12pt]{article}
\usepackage{amsmath,amssymb}
\usepackage{amsthm}
\usepackage{graphicx}
\usepackage{enumerate}
\usepackage{float}%  for figures
\usepackage{pstricks}%		PSTricks package
\usepackage{pst-eps}%
\usepackage{pst-plot}%
\usepackage{pst-node}%
\newtheorem{theorem}{Theorem}[section]
\newtheorem{proposition}[theorem]{Proposition}
\theoremstyle{definition}
\newtheorem{definition}[theorem]{Definition}
\theoremstyle{remark}
\newtheorem{remark}[theorem]{Remark}
\theoremstyle{example}
\newtheorem{example}[theorem]{Example}
\newtheorem{observation}[theorem]{Observation}
\theoremstyle{lemma}
\newtheorem{lemma}[theorem]{Lemma}
\theoremstyle{setting}

\theoremstyle{corollary}
\newtheorem{cor}[theorem]{Corollary}
\theoremstyle{problem}

\newcommand{\en}[2]{M^{#1}(#2)}
\newcommand{\ensub}[3]{M^{#1}_{#3}(#2)}
\numberwithin{equation}{section}
\numberwithin{figure}{section}
\makeatletter
\renewcommand{\@makecaption}[2]{%
  \vskip\abovecaptionskip
  \sbox\@tempboxa{#1 #2}%
  \ifdim \wd\@tempboxa >\hsize
    #1: #2\par
  \else
    \global \@minipagefalse
    \hb@xt@\hsize{\hfil\box\@tempboxa\hfil}%
  \fi
  \vskip\belowcaptionskip}
\makeatother
\usepackage[breaklinks]{hyperref}
\begin{document}

\title{Restrictions of generalized Verma modules\\
to symmetric pairs}
\author{Toshiyuki Kobayashi\footnote{Partially supported by
Institut des Hautes \'{E}tudes Scientifiques, France and
        Grant-in-Aid for Scientific Research (B) (22340026), Japan
        Society for the Promotion of Science}}
\date{}

\maketitle

\begin{abstract}
We initiate a new line of investigation on branching problems for
generalized Verma modules with respect to reductive symmetric pairs
$(\mathfrak{g},\mathfrak{g}')$.
In general, 
Verma modules may not contain any simple module when
restricted to a reductive subalgebra.
In this article
 we give a necessary and sufficient condition on the triple
$(\mathfrak{g},\mathfrak{g}',\mathfrak{p})$
such that the restriction $X|_{\mathfrak{g}'}$ always contains simple
$\mathfrak{g}'$-modules for any ${\mathfrak{g}}$-module $X$ lying in 
the parabolic BGG category $\mathcal{O}^{\mathfrak{p}}$ attached to a
parabolic subalgebra $\mathfrak{p}$ of $\mathfrak{g}$. 
Formulas are derived 
 for the Gelfand--Kirillov dimension
 of any simple module occurring in a simple generalized Verma module.   
We then prove that the restriction $X|_{\mathfrak{g}'}$ is generically
multiplicity-free for any $\mathfrak{p}$ and
any $X \in \mathcal{O}^{\mathfrak{p}}$
if and only if $(\mathfrak{g},\mathfrak{g}')$ is isomorphic to
$(A_n,A_{n-1})$, $(B_n,D_n)$, or $(D_{n+1},B_n)$.
Explicit branching laws are also presented.
\end{abstract}



\setcounter{tocdepth}{1}
\tableofcontents

\addtocounter{section}{0}

\section{Program}
\label{sec:1}

Branching problems in representation theory ask how irreducible
modules decompose when restricted to subalgebras.  
In the context of the Bernstein--Gelfand--Gelfand category $\mathcal{O}$ of a semisimple
Lie algebra $\mathfrak{g}$, 
 branching problems are seemingly simple,
 however, 
 it turns out that the restrictions behave wildly in general.   
For instance, 
  the restrictions $X|_{\mathfrak{g}_1'}$
 and $X|_{\mathfrak{g}_2'}$ of a $\mathfrak{g}$-module $X$ lying in
$\mathcal{O}$ may be completely different 
 even when two reductive subalgebras ${\mathfrak {g}}_1'$ and ${\mathfrak {g}}_2'$
 are conjugate to each other as the following observations indicate
 (see Examples \ref{ex:Heisen1}, \ref{ex:GKCA} for more details):

\begin{observation}
\label{ex:dc}
The restriction
$X|_{\mathfrak {g}_1'}$
 does NOT contain any simple $\mathfrak {g}_1'$-module,
whereas
 $X|_{\mathfrak {g}_2'}$
 decomposes into an algebraic direct sum 
 of simple ${\mathfrak {g}}_2'$-modules.  
\end{observation}

\begin{observation}
\label{ex:GKdim}
The Gelfand--Kirillov dimension 
 of any simple $
{\mathfrak {g}}_1'
$-module occurring in $X|_{\mathfrak {g}_1'}$
 is larger than that of any simple ${\mathfrak {g}}_2'$-module 
 in $X|_{\mathfrak {g}_2'}$.  
\end{observation}

The understanding of such phenomena requires a precise formulation of
 branching problems.
Among others we begin by asking what is
 a `well-posed' framework of branching problems for the restriction
$X|_{\mathfrak{g}'}$ where $\mathfrak{g}'$ is a (generalized)
 reductive subalgebra of $\mathfrak{g}$
and $X$ lies in the category $\mathcal{O}$:
\begin{description}
\item[\textbf{Problem A.}]
When does the restriction $X|_{\mathfrak{g}'}$
 contain a simple ${\mathfrak{g}}'$-module?
\end{description}
Further,
 we raise the following problems
 when $X|_{\mathfrak{g}'}$ contains simple $\mathfrak{g}'$-modules.  
\begin{description}
\item[\textbf{Problem B.}]
Find the \lq{size}\rq\
 of simple ${\mathfrak{g}}'$-modules
occurring in $X|_{{\mathfrak{g}}'}$.  
\item[{\textbf{Problem C.}}]
Estimate multiplicities
 of simple ${\mathfrak {g}}'$-modules occurring 
 in  $X|_{{\mathfrak{g}}'}$.  
\item[{\textbf{Problem D.}}]
Find branching laws,
in particular, 
 for multiplicity-free cases.
\end{description}

Let us explain briefly our main results.  
We write $\mathfrak{B}$ for the full flag variety of $\mathfrak{g}$,
and $\mathfrak{G}'$ for the set of 
conjugacy classes of $\mathfrak{g}'$ under the group
$G := \operatorname{Int}(\mathfrak{g})$ of inner automorphisms. 
Then the `framework' of the restriction
$X|_{\mathfrak{h}}$ for $X \in \mathcal{O}$ and
$\mathfrak{h} \in \mathfrak{G}'$ is described by means of the quotient space
$G \backslash
 (\mathfrak{B} \times \mathfrak{G}')$
under the diagonal action of $G$.
More generally,
we formulate a proper framework to discuss Problems A to D in Theorem \ref{thm:orbit} in the
parabolic BGG category
$\mathcal{O}^{\mathfrak{p}}$
 (see Subsection \ref{subsec:2.1})
 for an arbitrary parabolic subalgebra ${\mathfrak {p}}$ of ${\mathfrak{g}}$.
After discussing basic results in this framework in the generality
that $\mathfrak{g}'$ is an arbitrary reductive subalgebra in
$\mathfrak{g}$, 
we highlight the case
 where $({\mathfrak {g}}, {\mathfrak {g}}')$ is a symmetric pair to
get finer results, keeping differential geometric applications in mind.
It includes the `group case' 
$
({\mathfrak {g}}_1 \oplus {\mathfrak {g}}_1,\operatorname{diag}({\mathfrak {g}}_1))
$
as a special example, 
for which the branching laws 
 describe the decomposition of the tensor product of two
representations (e.g.\ \textit{fusion rules}). 
For symmetric pairs $(\mathfrak{g},\mathfrak{g}')$,
the cardinality of $G$-orbits on
$\mathfrak{B} \times \mathfrak{G}'$ is finite, and
 we give a complete answer to Problem A in the category ${\mathcal{O}}^{\mathfrak {p}}$
 in terms of the finite set
$G \backslash (\mathfrak{P} \times
 \mathfrak{G}')$.
Namely,
 we prove in Theorem \ref{thm:decoc} that the restriction 
$X|_{\mathfrak{g}'}$ contains simple $\mathfrak{g}'$-modules 
for any $X \in \mathcal{O}^{\mathfrak{p}}$
if and only if $(\mathfrak{p},\mathfrak{g}')$ lies in a closed
$G$-orbit on
$\mathfrak{P} \times \mathfrak{G}'$.

Turning to Problem B, 
 we make use of the associated varieties (see e.g.\ \cite{Ja,V91}) as
 a coarse measure
 of the \lq{size}\rq\ of ${\mathfrak {g}}'$-modules.  
We see that
 the associated variety $\mathcal{V}_{\mathfrak{g}'}(Y)$
of a simple ${\mathfrak {g}}'$-module $Y$ occurring in the restriction
$X|_{\mathfrak{g}'}$ is independent of $Y$ if $X$ is a simple
$\mathfrak{g}$-module.  
The formulas of
$\mathcal{V}_{\mathfrak{g}'}(Y)$ 
and its dimension (\textit{Gelfand--Kirillov dimension}) are derived
in Theorem \ref{thm:AV}.

Concerning Problem C,
it is notorious in the category of unitary representations of real
reductive groups that the multiplicities in the branching laws
may be infinite
 when restricted to symmetric pairs,
 see \cite{Ko2}.  
In contrast, 
 we prove in Theorem \ref{thm:fm} 
 that multiplicities are always finite
in the branching laws with respect to symmetric pairs
 in the category ${\mathcal{O}}$.

Particularly interesting branching laws are multiplicity-free cases
 where any simple ${\mathfrak {g}}'$-module occurs in the restriction $X |_{\mathfrak{g}'}$ 
at most once.  
We  give two general multiplicity-free theorems
 with respect to symmetric pairs $(\mathfrak{g},\mathfrak{g}')$
 in the parabolic category ${\mathcal{O}}^{\mathfrak{p}}$:
\begin{enumerate}
\item[1)] 
${\mathfrak {p}}$ special, 
 $({\mathfrak {g}},{\mathfrak {g}}')$ general
  (Theorem \ref{thm:PAbel}),
\item[2)] 
${\mathfrak {p}}$ general, 
 $({\mathfrak {g}},{\mathfrak {g}}')$ special (Theorem \ref{thm:mfall}),
\end{enumerate}
and then find  branching laws corresponding to closed orbits in
$G \backslash(\mathfrak{P}\times\mathfrak{G}')$.

This is the first article of our project on a systematic construction
of equivariant differential operators in parabolic geometry.
In subsequent papers,
Theorem \ref{thm:decoc} (a solution to Theorem A) plays a foundational
role in dealing with 
\begin{itemize}
\item
a construction of conformally equivariant differential operators in
parabolic geometry,
\item
a generalization of the
 Rankin--Cohen bracket operators.
\end{itemize}
Actual calculations are carried out by using algebraic branching formulas
(Theorem \ref{thm:mfbr1}) together with an analytic machinery that we
call the `$F$-method' in \cite{koss}.

In Section \ref{subsec:PAbel} we have studied parabolic subalgebras
$\mathfrak{p}$ with abelian nilpotent radical.
The case of parabolic subalgebras $\mathfrak{p}$ with Heisenberg
nilpotent radical (e.g.\ Example \ref{ex:Heisen1}) may be thought of
as a generalization of Section \ref{subsec:PAbel}.
Using Theorems \ref{thm:br} and \ref{thm:decoc},
we can give a complete classification of the triples
$(\mathfrak{g},\mathfrak{p},\mathfrak{g}^\tau)$ 
and the closed orbits in $G^\tau\backslash G/P$ (see the framework of
Theorem \ref{thm:orbit})
with discrete
decomposable and multiplicity-free
branching laws.
The calculation is more involved,
and will be reported in another paper.

Partial results of this article were presented
 at the conference
 in honor of Vinberg's 70th birthday
 at Bielefeld in Germany in 2007
 and a series of lectures
at the Winter School on Geometry and Physics in Cech Republic in 2010.  
The author is grateful to the organizers, in particular,
Professors Abels and Sou\v{c}ek,
 for their warm hospitality.  

{\bf{Notation:}}\enspace
${\mathbb{N}}=\{0,1,2, \cdots\}$, 
${\mathbb{N}}_+=\{1,2,3, \cdots\}$.

\section{Branching problem of Verma modules}
\label{sec:2}

In general, Verma modules 
may not contain any simple $\mathfrak{g}'$-module
when restricted to a reductive subalgebra $\mathfrak{g}'$.
In this section,
we use the geometry of the double coset space
$N_G(\mathfrak{g}') \backslash G / P$
and clarify the problem in Theorem \ref{thm:orbit},
which will then serve as a foundational setting of branching
problems for the category
$\mathcal{O}^{\mathfrak{p}}$ in Theorem \ref{thm:decoc}.

\subsection{Generalized Verma modules}
\label{subsec:2.1}

We begin with a quick review of the (parabolic) BGG category
$\mathcal{O}^{\mathfrak{p}}$
 and fix some notation.

Let ${\mathfrak {g}}$ be a semisimple Lie 
 algebra over $\mathbb C$, 
 and ${\mathfrak {j}}$ a Cartan subalgebra.  
We write $\Delta \equiv \Delta({\mathfrak {g}},
{\mathfrak {j}})$
 for the root system, 
 ${\mathfrak {g}}_{\alpha}$
 ($\alpha \in\Delta$) for the root space, 
 and $\alpha^{\vee}$ for the coroot. 
We fix a positive system $\Delta^+$, 
 and define a Borel subalgebra
 ${\mathfrak {b}}={\mathfrak {j}}
+{\mathfrak {n}}$
 with nilradical  ${\mathfrak {n}}
:= \oplus_{\alpha\in \Delta^+}
{\mathfrak {g}}_{\alpha}$.   
The BGG category ${\mathcal {O}}$ is defined
 to be the full subcategory 
 of ${\mathfrak {g}}$-modules
 whose objects are finitely generated ${\mathfrak {g}}$-modules $X$
 such that 
 $X$ are ${\mathfrak {j}}$-semisimple 
 and locally ${\mathfrak {n}}$-finite \cite{BGG}.

Let $\mathfrak{p}$ be a standard parabolic subalgebra,
and ${\mathfrak{p}} = {\mathfrak{l}} + {\mathfrak{u}}$
 its Levi decomposition 
 with ${\mathfrak{j}} \subset {\mathfrak{l}}$.
We set 
$\Delta^+({\mathfrak{l}}):=\Delta^+ \cap \Delta({\mathfrak{l}}, {\mathfrak{j}})$,
and define
$
{\mathfrak {n}}_-({\mathfrak {l}})
  :=\oplus_{\alpha\in \Delta^+ ({\mathfrak {l}})}
  {\mathfrak {g}}_{-\alpha}$.
The parabolic BGG category 
${\mathcal{O}}^{\mathfrak {p}}$
 is defined to be the full subcategory 
 of ${\mathcal{O}}$
 whose objects $X$ are locally 
 ${\mathfrak {n}}_-({\mathfrak {l}})$-finite.  
Then ${\mathcal{O}}^{\mathfrak{p}}$ is closed
 under submodules, quotients, and tensor products
 with finite dimensional representations.

The set of $\lambda$ for 
 which $\lambda|_{{\mathfrak{j}} \cap [{\mathfrak {l}},{\mathfrak {l}}]}$
 is dominant integral is denoted by 
\[
  \Lambda^+({\mathfrak {l}})
:=
 \{\lambda  \in {\mathfrak {j}}^*
  :
   \langle \lambda,\alpha^\vee\rangle
   \in {\mathbb{N}}
   \text{ for all }\alpha\in \Delta^+({\mathfrak {l}})\}. 
\]
We write $F_\lambda$ for the finite dimensional simple ${\mathfrak{l}}$-module
 with highest weight $\lambda$,
 inflate $F_\lambda$ to a ${\mathfrak {p}}$-module
 via the projection ${\mathfrak {p}}\to {\mathfrak {p}}/{\mathfrak {u}}\simeq {\mathfrak {l}}$, 
 and define the generalized Verma module by 
\begin{equation}
\label{eqn:NGpF}
  \ensub {\mathfrak{g}} \lambda {\mathfrak {p}}
  \equiv \ensub{\mathfrak {g}}{F_{\lambda}}{\mathfrak {p}}
:= U({\mathfrak {g}}) \otimes_{U({\mathfrak{p}})} F_{\lambda}.  
\end{equation}
Then $\ensub{\mathfrak {g}}{\lambda} {\mathfrak {p}}\in
{\mathcal{O}}^{\mathfrak {p}}$,
and any simple object in $\mathcal{O}^{\mathfrak {p}}$ is the quotient
of some $\ensub{\mathfrak {g}}{\lambda} {\mathfrak {p}}$.
We say $\ensub{\mathfrak {g}}{\lambda}{\mathfrak {p}}$ is of {\it{scalar type}}
 if $F_{\lambda}$ is one-dimensional, or equivalently,
 if $\langle \lambda, \alpha^\vee\rangle =0$
 for all $\alpha \in \Delta({\mathfrak {l}})$.

Let $\rho$ be half the sum of positive roots.  
If $\lambda \in \Lambda^+({\mathfrak {l}})$
 satisfies 
\begin{equation}
\label{eqn:anti}
   \langle \lambda + \rho , \beta^\vee\rangle \not\in {\mathbb{N}}_+
  \text{ for all } \beta \in \Delta^+ \setminus \Delta({\mathfrak {l}}), 
\end{equation}
then $\ensub{\mathfrak {g}}{\lambda}{\mathfrak {p}}$ is simple, 
 see \cite{BD}.

For ${\mathfrak{p}}={\mathfrak{b}}$, 
 we simply write $\en{{\mathfrak {g}}}{\lambda}$
 for $\ensub{{\mathfrak {g}}}{\lambda}{{\mathfrak {b}}}$.
We note that ${\mathcal{O}}^{\mathfrak{b}}={\mathcal{O}}$
 by definition.

\subsection{Framework of branching problems}
\label{subsec:form}

Let ${\mathfrak{g}}'$ be a subalgebra of ${\mathfrak{g}}$, 
 and ${\mathfrak{p}}$ a parabolic subalgebra of ${\mathfrak{g}}$.  
We denote by ${\mathfrak{G}}'$ and $\mathfrak{P}$ the set of conjugacy classes
 of ${\mathfrak{g}}'$
and $\mathfrak{p}$, respectively.
Let
 $P$ be the parabolic subgroup of 
$G = \operatorname{Int}(\mathfrak{g})$
with Lie algebra
 ${\mathfrak {p}}$,
and define the normalizer of ${\mathfrak{g}}'$ as  
\[
  N_{G}(\mathfrak {g}'):= 
 \{g \in G: \operatorname{Ad} (g) {\mathfrak {g}}' = {\mathfrak {g}}'\}.  
\]
Then we have natural bijections:
$G/P \simeq \mathfrak{P}$,
$G/N_G(\mathfrak{g'}) \simeq \mathfrak{G}'$,
and hence
\begin{equation}\label{eqn:ngp}
G \backslash (\mathfrak{P} \times \mathfrak{G}')
\simeq
N_G(\mathfrak{g'}) \backslash \mathfrak{P}
\simeq
\mathfrak{G}' / P
\simeq
N_G(\mathfrak{g}') \backslash G / P.
\end{equation}
Here, we let $G$ act diagonally on
$\mathfrak{P} \times \mathfrak{G}'$
in the left-hand side of \eqref{eqn:ngp}.

Let $S$ be the set of complete representatives of the double coset
$N_G(\mathfrak{g}') \backslash G / P$, 
 and we write ${\mathfrak{g}}_s':=\operatorname{Ad}(s)^{-1}{\mathfrak{g}}'$
 for  $s \in S$.  
Then the branching problem
 for ${\mathcal {O}}^{\mathfrak {p}}$
 with respect to a subalgebra belonging to ${\mathfrak {G}}'$
 is \lq{classified}\rq\  
 by the double coset $N_G ({\mathfrak {g}}') \backslash G/P$ in the following sense:
\begin{theorem}
\label{thm:orbit}
For any $X \in \mathcal{O}^{\mathfrak{p}}$ and 
any $\mathfrak{h} \in \mathfrak{G}'$,
there exists
 $s \in S$ such that
$X|_{\mathfrak{h}} \simeq {\widetilde{X}}|_{\mathfrak{g}_s'}$
 for some ${\widetilde{X}} \in {\mathcal{O}}^{\mathfrak{p}}$
via a Lie algebra isomorphism between $\mathfrak{h}$ and $\mathfrak{g}_s'$.
\end{theorem}

\begin{proof}
[Proof of Theorem \ref{thm:orbit}]
Given $\mathfrak{h} \in \mathfrak{G}'$,
we take $s \in S$ and $q \in P$ such that
$\operatorname{Ad}((sq)^{-1})\mathfrak{g}' = \mathfrak{h}$.
Clearly,
we have a Lie algebra isomorphism
$\operatorname{Ad}(q^{-1}): \mathfrak{g}_s' \overset{\sim}{\to} \mathfrak{h}$.

For $X \in {\mathcal{O}}^{\mathfrak{p}}$,  
 we define a new ${\mathfrak{g}}$-module structure on $X$
 by 
\[     
  Z \underset q \cdot v := (\operatorname{Ad}(q)^{-1} Z) \cdot v
 \quad \text{for } Z \in {\mathfrak {g}}, v \in X.  
\]

Since $P$ normalizes $\mathfrak{p}$,
this new module, 
 to be denoted by $\widetilde X$, 
 lies in ${\mathcal {O}}^{\mathfrak {p}}$.  
Then, for any Lie subalgebra $\mathfrak{v}$ of $\mathfrak{g}$,
the restriction ${\widetilde{X}}|_{\mathfrak{v}}$ is isomorphic to
the restriction
$X|_{\operatorname{Ad}(q)^{-1}\mathfrak{v}}$
via the Lie algebra isomorphism
$\mathfrak{v} \simeq \operatorname{Ad}(q)^{-1}\mathfrak{v}$.
Applying this to 
$\mathfrak{v} := \mathfrak{g}_s'$,
we get the following isomorphism:
\[
X|_{\mathfrak{h}}
=
X|_{\operatorname{Ad}(q)^{-1}\operatorname{Ad}(s)^{-1}\mathfrak{g}'}
\simeq
{\widetilde{X}}|_{\mathfrak{g}_s'}
\]
via the Lie algebra isomorphism 
$\operatorname{Ad}(q) : {\mathfrak{h}} \overset{\sim}{\rightarrow} {\mathfrak{g}}_s'$.  
Theorem \ref{thm:orbit} is thus proved.
\end{proof}

\begin{remark}
If $(\mathfrak{g},\mathfrak{g}')$ is a semisimple symmetric pair
(see Subsection \ref{subsec:2.4}),
then $S$ is a finite set (Matsuki \cite{M79}).
\end{remark}

\section{Discretely decomposable branching laws}
\label{sec:deco}

In this section, 
 we bring the concept of 
\lq{discretely decomposable restrictions}\rq\
to the branching problem
 for the BGG category ${\mathcal{O}}^{\mathfrak{p}}$,
 and prove that the restriction $X|_{\mathfrak{g}'}$ 
 contains simple ${\mathfrak{g}}'$-modules
 for $X \in {\mathcal{O}}^{\mathfrak{p}}$
 if ${\mathfrak{p}}$ lies in a closed $G'$-orbit 
 on the  generalized flag variety ${\mathfrak{P}}$.
In particular, 
 it is the case
 if ${\mathfrak{p}}$ is ${\mathfrak {g}}'$-compatible 
 (Definition \ref{def:cp}).
Under this assumption the character identities are derived
 for the restriction $X|_{\mathfrak{g}'}$
 (Theorem \ref{thm:br}).

\subsection{Discretely decomposable modules
 ${\mathcal{O}}$}
\label{subsec:2.2}

Suppose that $\mathfrak{g}$ is a reductive Lie algebra.
\begin{definition}
\label{def:discdeco}
We say a $\mathfrak{g}$-module $X$ is \textit{discretely decomposable} if there is an
increasing filtration $\{X_m\}$ of $\mathfrak{g}$-submodules of finite length 
such that 
$X = \bigcup_{m=0}^\infty X_m$.  
Further,
we say $X$ is {\it{discretely decomposable
 in the category}} ${\mathcal {O}}^{\mathfrak{p}}$
 if all $X_m$ can be taken from ${\mathcal{O}}^{\mathfrak{p}}$.  

\end{definition}

Here are obvious examples:
\begin{example}
\label{ex:discdeco}
{\rm{1)}}\enspace
Any $\mathfrak{g}$-module of finite length is 
discretely decomposable.
\par\noindent
{\rm{2)}}\enspace
{\rm{(completely reducible case).}}
\enspace
An algebraic direct sum of countably many simple ${\mathfrak {g}}$-modules
 is discretely decomposable.  
\end{example}

\begin{remark}
\label{rem:discdeco}
The concept of discretely decomposable ${\mathfrak{g}}$-modules 
 was originally introduced 
 in the context of $({\mathfrak{g}},K)$-modules
 in {\cite[Definition 1.1]{Ko3}}
 as an algebraic analogue
 of  unitary representations whose irreducible decompositions have no
\lq{continuous spectrum}\rq.  
Then the main issue of \cite{Ko2,Ko3} was to find a criterion for the discrete
 decomposability
 of the restriction of $({\mathfrak{g}}, K)$-modules.
We note 
 that discrete decomposability in the generality of Definition
 \ref{def:discdeco} does not imply complete reducibility. 
\end{remark}
\medskip
Suppose ${\mathfrak{g}}'$ is a reductive subalgebra,
 and ${\mathfrak{p}}'$  its parabolic subalgebra.

\begin{lemma}
\label{lem:YX}
Let $X$ be a simple ${\mathfrak {g}}$-module.   
Then the restriction $X|_{\mathfrak {g}'}$
 is discretely decomposable in the category ${\mathcal{O}}^{\mathfrak{p}'}$
 if and only if there exists a ${\mathfrak {g}}'$-module $Y \in {\mathcal {O}}^{\mathfrak{p}'}$
 such that 
$
     \operatorname{Hom}_{\mathfrak {g}'} (Y,X|_{\mathfrak {g}'}) \ne \{0\}.  
$
In this case,
any subquotient occurring in the $\mathfrak{g}'$-module $X|_{\mathfrak{g}'}$
lies in $\mathcal{O}^{\mathfrak{p}'}$.
\end{lemma}

\begin{proof}
Suppose
$\operatorname{Hom}_{\mathfrak{g}'}(Y,X|_{\mathfrak{g}'}) \ne 0$
for some $Y \in \mathcal{O}^{\mathfrak{p}'}$.
Taking the subquotient of $Y$ if necessary,
we may assume $Y$ is a simple ${\mathfrak {g}}'$-module.  
Let $\iota: Y \to X$ be an injective $\mathfrak{g}'$-homomorphism.
For $m \in {\mathbb{N}}$, 
 we denote by $Y_m$ the image of the following
$\mathfrak{g}'$-homomorphism:
\[
\underbrace{\mathfrak{g} \otimes \cdots \otimes \mathfrak{g}}_m
\otimes Y \to X,
\quad
(H_1 \otimes \cdots \otimes H_m) \otimes v \mapsto
H_1 \cdots H_m \, \iota(v).
\]
Then $X = \bigcup_{m=0}^\infty Y_m$ because $X$ is simple.
Moreover $Y_m \in \mathcal{O}^{\mathfrak{p}'}$
 because ${\mathcal{O}}^{\mathfrak{p}'}$ is closed under quotients
 and tensor products with finite dimensional representations.
Hence,
 the restriction $X|_{\mathfrak{g}'}$ is discretely decomposable in
 $\mathcal{O}^{\mathfrak{p}'}$.
Conversely, the `only if' part is obvious because
 ${\mathcal{O}}^{\mathfrak {p}'}$ is closed
under submodules.
Finally,
 any subquotient of $Y_m$ lies in ${\mathcal{O}}^{\mathfrak{p}'}$,
 whence the last statement.    
Thus Lemma \ref{lem:YX} is proved.
\end{proof}

\subsection{Discretely decomposable restrictions for $\mathcal{O}^{\mathfrak{p}}$}
\label{subsec:closed}
Let $G=\operatorname{Int}({\mathfrak{g}})$, 
 $P$ the parabolic subgroup
 of $G$ with Lie algebra ${\mathfrak{p}}$ as before, 
 and $G'$ a reductive subgroup with Lie algebra ${\mathfrak {g}}'$.  
We ask when the restriction $X|_{\mathfrak{g}'}$ of 
 $X \in {\mathcal{O}}^{\mathfrak{p}}$
is discretely decomposable in the sense of Definition \ref{def:discdeco}.
\begin{proposition}
\label{prop:closed}
 If $G'P$ is closed in $G$
 then the restriction $X|_{{\mathfrak {g}}'}$ is discretely decomposable 
 for any simple ${\mathfrak {g}}$-module $X$ in ${\mathcal{O}}^{\mathfrak {p}}$.  
\end{proposition}

\begin{proof}
We set $P' := G' \cap P$.    
Suppose $G' P$ is closed in $G$.
Then $G'/P'$ is closed in the generalized flag variety $G/P$, 
 and hence is compact.  
Therefore,
 the Lie algebra ${\mathfrak {p}}' := {\mathfrak {g}}' \cap {\mathfrak {p}}$ of $P'$
 must be a parabolic subalgebra of ${\mathfrak {g}}'$.

Let $X$ be a simple object in ${\mathcal{O}}^{\mathfrak {p}}$.  
Then $X$ is obtained as the quotient of some generalized Verma module,
 that is,
 there exists $\lambda \in \Lambda^+ ({\mathfrak{l}})$
 such that the composition map 
\begin{equation*}
     F_{\lambda}
    \hookrightarrow 
    U({\mathfrak {g}}) \otimes_{U({\mathfrak {p}})}
     F_{\lambda}
    \twoheadrightarrow 
    X
\end{equation*}
 is non-trivial.  
Therefore,
 we get a non-zero ${\mathfrak{g}}'$-homomorphism
\begin{equation}
\label{eqn:YXlmd}
U({\mathfrak {g}}') \otimes_{U({\mathfrak {p}}')}
     (F_{\lambda}|_{\mathfrak {p}'})
  \to X.  
\end{equation}
Since the $\mathfrak{g}'$-module 
$U({\mathfrak {g}}') \otimes_{U({\mathfrak {p}}')}
     (F_{\lambda}|_{\mathfrak {p}'})$ lies in ${\mathcal{O}}^{\mathfrak{p}'}$,  
 the restriction $X|_{\mathfrak {g}'}$
 is discretely decomposable
 in the category ${\mathcal{O}}^{\mathfrak {p}'}$
 owing to Lemma \ref{lem:YX}.  
\end{proof}
The converse statement of Proposition \ref{prop:closed}
 will be proved in Theorem \ref{thm:decoc}
 under the assumption
 that $({\mathfrak{g}}, {\mathfrak{g}}')$ is a semisimple symmetric pair.

The assumption of Proposition \ref{prop:closed}
 fits well into the framework of Theorem \ref{thm:orbit}.  
To see this,
 we make the following (easy) observation:
\begin{lemma}
\label{lem:HPclosed}
Retain the notation as in Subsection \ref{subsec:form}.
Then the following conditions
 on the triple $({\mathfrak{g}}, {\mathfrak{g}}', {\mathfrak{p}})$
 are equivalent:
\begin{enumerate}
\item[{\rm{(i)}}]
The $G'$-orbit through ${\mathfrak{p}} \in {\mathfrak{P}}$ is closed.   
\item[{\rm{(ii)}}]
$G'P$ is closed in $G$.  
\end{enumerate}
\end{lemma}
Clearly these conditions are invariant
 under the conjugation of $\mathfrak{p}$ by an element of the group
$N_G({\mathfrak{g}}')$, 
 and hence they are determined by the equivalence classes
 in $N_G({\mathfrak{g}}')\backslash G/P \simeq G \backslash ({\mathfrak{P}} \times {\mathfrak{G}}')$
 (see \eqref{eqn:ngp})
 containing $({\mathfrak{p}}, {\mathfrak{g}}')\in {\mathfrak{P}} \times {\mathfrak{G}}'$.

\subsection{$\protect{\mathfrak{g}}'$-compatible parabolic subalgebra $\protect\mathfrak {p}$}
\label{subsec:cp}

This subsection discusses a sufficient condition
 for the closedness of $G'P$ in $G$.

A semisimple element $H \in {\mathfrak{g}}$
is said to be {\it{hyperbolic}}
 if the eigenvalues of $\operatorname{ad}(H)$
 are all real.  
For a hyperbolic element $H$,
 we define the subalgebras 
\[
 {\mathfrak{u}}_+ \equiv {\mathfrak{u}}_+(H),
 \quad
 {\mathfrak{l}}\equiv {\mathfrak{l}}(H), 
 \text{ and }
  {\mathfrak{u}}_-\equiv {\mathfrak{u}}_-(H)
\]
 as the sum of the eigenspaces
 with positive,
 zero, 
 and negative eigenvalues,
 respectively.  
Then 
\begin{equation}
\label{eqn:pH}
  {\mathfrak{p}}(H) :=
  {\mathfrak{l}}(H)
  +{\mathfrak{u}}_+(H)
\end{equation}
 is a Levi decomposition of a parabolic subalgebra
 of ${\mathfrak{g}}$. 

Let ${\mathfrak{g}}'$ be a reductive subalgebra
 of ${\mathfrak {g}}$,
 and ${\mathfrak{p}}$ a parabolic subalgebra
 of ${\mathfrak{g}}$. 

\begin{definition}
\label{def:cp}
We say ${\mathfrak {p}}$ is ${\mathfrak{g}}'$-\textit{compatible}
 if there exists a hyperbolic element
 $H$ in ${\mathfrak{g}}'$
 such that ${\mathfrak{p}}={\mathfrak{p}}(H)$.  
\end{definition}

If ${\mathfrak {p}}={\mathfrak{l}}+{\mathfrak{u}}_+$ is ${\mathfrak{g}}'$-compatible, 
then 
$
    {\mathfrak{p}}':={\mathfrak{p}} \cap {\mathfrak{g}}'
$
 becomes a parabolic subalgebra 
 of ${\mathfrak{g}}'$
 with Levi decomposition
\[
   {\mathfrak{p}}'={\mathfrak{l}}'+{\mathfrak{u}}_+'
                    :=({\mathfrak{l}} \cap {\mathfrak{g}}') + ({\mathfrak{u}}_+ \cap {\mathfrak{g}}').  
\]
Then, 
 using the notation of Subsection \ref{subsec:closed},
 we see that $G'/P'=G'/G' \cap P$ becomes a generalized flag variety,
 and therefore is closed in $G/P$.  
Hence,
 we get the following proposition from Proposition \ref{prop:closed}:
\begin{proposition}
\label{prop:cp}
If ${\mathfrak{p}}$ is ${\mathfrak{g}}'$-compatible,
 then $G'P$ is closed in $G$
 and the restriction $X|_{\mathfrak{g}'}$ is discretely decomposable for any $X \in {\mathcal{O}}^{{\mathfrak{p}}}$.  
\end{proposition}

The converse statement holds when $\mathfrak{p}$ is a Borel subalgebra
(Corollary \ref{cor:bdeco}),
but does not always hold for a general parabolic subalgebra.
Theorem \ref{thm:decoc} below shows that the following example gives a
counterexample to the converse statement of Proposition \ref{prop:cp}.

\begin{example}\label{ex:cp}
Let $\mathfrak{g} = \mathfrak{g}_1 \oplus \mathfrak{g}_1$,
and
$\mathfrak{g}' := \operatorname{diag}(\mathfrak{g}_1)
 \equiv \{(Z,Z): Z \in \mathfrak{g}_1\}$.
Then a parabolic subalgebra $\mathfrak{p}$ of\/ $\mathfrak{g}$ is 
$\mathfrak{g}'$-compatible if and only if\/ $\mathfrak{p}$ is of the
form $\mathfrak{p}_1 \oplus \mathfrak{p}_1$ for some parabolic
subalgebra $\mathfrak{p}_1$ of\/ $\mathfrak{g}_1$.

On the other hand,
 $G'P$ is closed in $G = G_1 \times G_1$ if and only if\/ $\mathfrak{p}$
is of the form $\mathfrak{p}_1 \oplus \mathfrak{p}_2$ for some
parabolic subalgebras $\mathfrak{p}_1$ and\/ $\mathfrak{p}_2$ containing
a common Borel subalgebra.
\end{example}

\subsection{Character identities}
\label{subsec:chid}

In this subsection,
 we prove the character identities
 of the restriction
 of generalized Verma modules
 to a reductive subalgebra $\mathfrak{g}'$ assuming that the
  parabolic subalgebras ${\mathfrak {p}}$ is $\mathfrak{g}'$-compatible.  

Let ${\mathfrak{p}}={\mathfrak{l}}+{\mathfrak{u}}_+$ be a ${\mathfrak{g}}'$-compatible parabolic subalgebra
 of ${\mathfrak{g}}$
 defined by a hyperbolic element $H \in {\mathfrak {g}}'$.  
We take a Cartan subalgebra ${\mathfrak {j}}'$ of ${\mathfrak {g}}'$
 such that $H \in {\mathfrak {j}}'$,
and extend it  to a Cartan subalgebra ${\mathfrak {j}}$ of ${\mathfrak {g}}$.  
Clearly,
 ${\mathfrak {j}} \subset {\mathfrak {l}}$
 and ${\mathfrak {j}}' \subset {\mathfrak {l}}'$.

We recall that $F_{\lambda}$ denotes the finite dimensional,
 simple module
 of ${\mathfrak {l}}$
 with highest weight $\lambda \in \Lambda^+({\mathfrak {l}})$.  
Likewise,
 let $F_{\delta}'$ denote that of ${\mathfrak {l}}'$
 for $\delta \in \Lambda^+({\mathfrak {l}'})$.  

Given a vector space $V$ we denote by
$S(V) = \bigoplus_{k=0}^\infty S^k(V)$
the symmetric tensor algebra over $V$.
We extend the adjoint
action of ${\mathfrak{l}}'$ on ${\mathfrak {u}}_-/{\mathfrak{u}}_- \cap {\mathfrak{g}}'$
 to $S({\mathfrak {u}}_-/{\mathfrak {u}}_- \cap {\mathfrak {g}}')$.  
We set 
\begin{equation}
\label{eqn:mab}
     m(\delta; \lambda):= \dim \operatorname{Hom}_{{\mathfrak {l}}'} 
                          (F_{\delta}', F_{\lambda}|_{\mathfrak{l}'} \otimes S({\mathfrak {u}}_-/{\mathfrak {u}}_- \cap {\mathfrak {g}}')).  
\end{equation}
\begin{theorem}
\label{thm:br}
Suppose that ${\mathfrak {p}}={\mathfrak{l}}+{\mathfrak{u}}_+$
 is a ${\mathfrak{g}}'$-compatible parabolic subalgebra of ${\mathfrak{g}}$, 
and $\lambda \in \Lambda^+({\mathfrak{l}})$.  

\par\noindent
{\rm{1)}}\enspace
$m(\delta;\lambda) < \infty$ for all $\delta \in \Lambda^+({\mathfrak{l}}')$.  
\par\noindent
{\rm{2)}}\enspace
In the Grothendieck group of $\mathcal{O}^{\mathfrak{p}'}$, 
 we have the following isomorphism:
\begin{equation}
\label{eqn:Mch}
  \ensub{\mathfrak{g}}{\lambda}{\mathfrak{p}}|_{{\mathfrak{g}}'} 
  \simeq
  \bigoplus_{\delta \in \Lambda^+({\mathfrak{l}}')}
  m(\delta;\lambda)
  \ensub{{\mathfrak{g}}'}{\delta}{\mathfrak{p}'}.  
\end{equation}
\end{theorem}
\begin{proof}
Let $H \in \mathfrak{g}'$ be the hyperbolic element defining the
parabolic subalgebra $\mathfrak{p}$.
We denote by
${\mathfrak{g}}''$  the orthogonal complementary subspace of ${\mathfrak{g}}'$
 in ${\mathfrak{g}}$
 with respect to the Killing form.  
Since $\operatorname{ad}(H)$ preserves the decomposition ${\mathfrak{g}}= {\mathfrak{g}}' \oplus {\mathfrak{g}}''$, 
the sum $\mathfrak{u}_-$ of negative eigenspaces of $\operatorname{ad}(H)$ decomposes as
\begin{equation}
\label{eqn:upm}
{\mathfrak{u}}_- = {\mathfrak{u}}_-' \oplus {\mathfrak{u}}_-''
:= ({\mathfrak{u}}_- \cap {\mathfrak{g}}')
\oplus ({\mathfrak{u}}_- \cap {\mathfrak{g}}'').
\end{equation}
This is a decomposition of ${\mathfrak{l}}'$-modules, 
 and hence, 
 we have an ${\mathfrak{l}}'$-module isomorphism
$
     S({\mathfrak{u}}_-'') \simeq S({\mathfrak{u}}_- / {\mathfrak{u}}_- \cap {\mathfrak{g}}')
$.

\par\noindent
{1)}\enspace
Let $a (>0)$ be the minimum of the eigenvalues of $-\operatorname{ad}(H)$ on ${\mathfrak {u}}_-''$.  
Since $H \in {\mathfrak{l}}'$, 
 we have  
\[
\operatorname{Hom}_{\mathfrak{l}'}(F_{\delta}', F_{\lambda} \otimes S^k({\mathfrak{u}}_-''))=0
\]
for all $k$ such that $k > \frac 1 a (\lambda(H) - \delta(H))$.  
In view of \eqref{eqn:mab}, 
we get  
$m(\delta;\lambda)< \infty$.  
\par\noindent
{2)}\enspace
The formal character of the generalized Verma module $\ensub {{\mathfrak{g}}}{\lambda}{\mathfrak{p}}$
 is given by 
\begin{equation}
\label{eqn:chM1}
  \operatorname{ch}(\ensub{{\mathfrak{g}}}{\lambda}{{\mathfrak{p}}})
=
  \operatorname{ch}(F_{\lambda})\prod_{\alpha \in \Delta({\mathfrak{u}}_-, {\mathfrak{j}})}(1-e^{\alpha})^{-1}.
\end{equation}
Let us prove that its restriction to ${\mathfrak {j}}'$
 equals the formal character of the right-hand side of \eqref{eqn:Mch}. 
For this,
we observe that
 $F_{\lambda} \otimes S({\mathfrak{u}}_-'')$ is a semisimple ${\mathfrak{l}}'$-module,
and therefore,
 it decomposes into the direct sum of simple ${\mathfrak{l}}'$-modules
$
     \bigoplus_{\delta \in \Lambda^+({\mathfrak{l}}')} m(\delta;\lambda)F_{\delta}',  
$
 where $m(\delta,\lambda)$ is defined in \eqref{eqn:mab}.  
Turning to their formal characters,
 we get 
\begin{equation}
\label{eqn:chM2}
   \operatorname{ch}(F_{\lambda})|_{{\mathfrak{j}}'}\prod_{\alpha \in \Delta({\mathfrak{u}}_-'', {\mathfrak{j}}')}(1-e^{\alpha})^{-1}
   =\sum_{\delta \in  \Lambda ^+({\mathfrak{l}}')} m(\delta;\lambda)\operatorname{ch}(F_{\delta}').  
\end{equation}
Writing the multiset $\Delta({\mathfrak{u}}_-, {\mathfrak{j}})|_{{\mathfrak{j}}'}$ as a disjoint union
$
   \Delta({\mathfrak{u}}_-'',{\mathfrak{j}}') \amalg \Delta({\mathfrak{u}}_-',{\mathfrak{j}}'), 
$
 we get from \eqref{eqn:chM1} and \eqref{eqn:chM2}
\begin{align*}
\operatorname{ch}(\ensub{\mathfrak{g}}{\lambda}{\mathfrak{p}}) |_{{\mathfrak{j}}'}
=&\operatorname{ch}(F_{\lambda})|_{{\mathfrak{j}}'} 
  \prod_{\alpha \in \Delta({\mathfrak{u}}_-'',{\mathfrak{j}}')}(1-e^{\alpha})^{-1}
  \prod_{\alpha \in \Delta({\mathfrak{u}}_-',{\mathfrak{j}}')}(1-e^{\alpha})^{-1}
\\
=&\sum_{\delta} m(\delta;\lambda) \operatorname{ch}(F_{\delta}')\prod_{\alpha \in \Delta({\mathfrak{u}}_-',{\mathfrak{j}}')}(1-e^{\alpha})^{-1}
\\
 =&\sum_{\delta} m(\delta;\lambda) \operatorname{ch}(\ensub{\mathfrak{g}'}{\delta}{\mathfrak{p}'}).  
\end{align*}
Hence 
 \eqref{eqn:Mch} holds in the Grothendieck group of ${\mathcal{O}}^{{\mathfrak{p}}'}$.  
\end{proof}

\subsection{Multiplicity-free restriction}
\label{subsec:mfbr}

Retain the setting of the previous subsection.  
In particular,
 we suppose that ${\mathfrak{p}}={\mathfrak{l}}+{\mathfrak{u}}_+$
 is a ${\mathfrak{g}}'$-compatible 
 parabolic  subalgebra of ${\mathfrak{g}}$.  
We will see in this subsection
that the character identity in Theorem \ref{thm:br}
 leads us to multiplicity-free branching laws
 for generalized Verma modules
 when the ${\mathfrak{l}}'$-module
 $S({\mathfrak{u}}_-/{\mathfrak{u}}_- \cap 
{\mathfrak{g}}')$ 
 is multiplicity-free.

\begin{definition}
\label{def:mfsp}
We say that a $\mathfrak{g}$-module $V$ is a \textit{multiplicity-free
space} if the induced $\mathfrak{g}$-module on the symmetric algebra
$S(V)$ is a multiplicity-free representation.
\end{definition}

Multiplicity-free spaces for reductive Lie algebras were classified by
V. Kac in the irreducible case, 
  and by Benson--Ratcliff 
 and Leahy
 independently in the reducible case (see \cite{BR}).

\medskip

The following Corollary is an immediate consequence of Theorem \ref{thm:br}:
\begin{cor}
\label{cor:br}
Assume that ${\mathfrak{u}}_-/{\mathfrak{u}}_- \cap {\mathfrak{g}}'$ is an ${\mathfrak{l}}'$-multiplicity-free space.  
We denote by $D$
 the support of simple ${\mathfrak{l}}'$-modules
 occurring in $S({\mathfrak{u}}_-/{\mathfrak{u}}_- \cap {\mathfrak{g}}')$, 
namely,
 $S({\mathfrak{u}}_-/{\mathfrak{u}}_- \cap {\mathfrak{g}}')
  \simeq 
   \bigoplus_{\delta\in D} F_{\delta}'$. 
Then any generalized Verma module $\ensub{{\mathfrak{g}}}{\lambda}{\mathfrak{p}}$
 of scalar type decomposes into a multiplicity-free sum
 of generalized Verma modules for ${\mathfrak{g}}'$
 in the Grothendieck group of ${\mathcal{O}}^{\mathfrak{p}'}$ as follows:
\begin{equation}
\label{eqn:Mmf}
     \ensub{\mathfrak{g}}{\lambda}{\mathfrak{p}}|_{\mathfrak{g}'}
     \simeq
     \bigoplus_{\delta \in D} \ensub{\mathfrak{g}'}{\lambda|_{\mathfrak{j}'}+ \delta}{\mathfrak{p}'}.     
\end{equation} 
\end{cor}

\begin{remark}
\label{rem:br}
For a \lq{generic}\rq\ $\lambda$,
 the formula \eqref{eqn:Mmf}
 becomes a multiplicity-free direct sum
 of simple ${\mathfrak{g}}'$-modules.  
For instance, 
there is no extension 
 among the modules $\ensub{\mathfrak{g}'}{\lambda|_{\mathfrak{j}'}+ \delta}{\mathfrak{p}'}$
 ($\delta \in D$)
 if they have distinct ${\mathfrak{Z}}({\mathfrak{g}}')$-infinitesimal characters
 (e.g.\ Theorems \ref{thm:AA}, \ref{thm:BD} and \ref{thm:DB})
 or if $\ensub{\mathfrak{g}}{\lambda}{\mathfrak{p}}$ has an invariant Hermitian inner product 
 with respect to a certain real form of ${\mathfrak {g}}'$
 (e.g.\ Theorem \ref{thm:PAbel}).  
See Section \ref{sec:mf} for details.  
\end{remark}

\section{Branching problems for symmetric pairs}

The decomposition of the tensor product of two representations
 is an example of branching laws with respect to
 a special case of symmetric pairs, namely,
 the pair
$\mathfrak{g}_1 \oplus \mathfrak{g}_1 \downarrow
 \operatorname{diag}(\mathfrak{g}_1)$.
In this section,
we discuss Problems A to D for semisimple symmetric pairs.

\subsection{Criterion for discretely decomposable restriction}
\label{subsec:2.4}

Let $\tau$ be an involutive automorphism 
 of a semisimple Lie algebra ${\mathfrak {g}}$, 
 and we denote the fixed point subalgebra by
\[
   {\mathfrak {g}}^{\tau}:=\{Z \in {\mathfrak {g}}: \tau Z=Z\}.  
\]
The pair $({\mathfrak {g}}, {\mathfrak {g}}^{\tau})$ is called a {\it{semisimple symmetric pair}}.  
Typical examples are the pairs
$(\mathfrak{g}_1 \oplus \mathfrak{g}_1, 
 \operatorname{diag}(\mathfrak{g}_1))$
($\mathfrak{g}_1$: semisimple Lie algebra),
$(\mathfrak{sl}_n, \mathfrak{so}_n)$,
and
$(\mathfrak{sl}_{p+q}, \mathfrak{s}(\mathfrak{gl}_p + \mathfrak{gl}_q))$.

We lift $\tau$ to an automorphism
 of the group $G=\operatorname{Int}(\mathfrak{g})$ of inner automorphisms, 
 and set $G^{\tau}:=\{g \in G:\tau g =g\}$.
Then $G^{\tau}$ is a reductive subgroup of $G$ with Lie algebra ${\mathfrak{g}}^{\tau}$.

Let $\mathfrak{p}$ be a parabolic subalgebra of $\mathfrak{g}$,
and $X$ a $\mathfrak{g}$-module lying in $\mathcal{O}^{\mathfrak{p}}$.
Problem A asks when
the restriction $X|_{\mathfrak{g}^{\tau}}$
 contains simple ${\mathfrak{g}}^{\tau}$-modules.
We give its necessary and sufficient condition by the geometry
of the generalized flag variety $G/P$ associated to the parabolic
 subalgebra $\mathfrak{p}$:

\begin{theorem}
\label{thm:decoc}
Let\/ $\mathfrak{g}$ be a complex semisimple Lie algebra,
$\tau$ an involutive automorphism of\/ $\mathfrak{g}$,
and\/ $\mathfrak{p}$ a parabolic subalgebra.
Then the following three conditions 
 on the triple $({\mathfrak {g}}, \mathfrak{g}^{\tau}, {\mathfrak{p}})$
 are equivalent:

\begin{enumerate}
\item[{\rm{(i)}}]
For any simple ${\mathfrak {g}}$-module $X$ in ${\mathcal{O}}^{\mathfrak {p}}$, 
 the restriction $X|_{{\mathfrak {g}}^{\tau}}$ contains 
 at least one simple ${\mathfrak {g}}^{\tau}$-module.  
\item[{\rm{(ii)}}]
For any simple ${\mathfrak {g}}$-module $X$ in ${\mathcal{O}}^{\mathfrak {p}}$, 
 the restriction $X|_{\mathfrak {g}^{\tau}}$ is discretely decomposable
 as a ${\mathfrak {g}}^{\tau}$-module
 in the sense of Definition \ref{def:discdeco}.  
\item[{\rm{(iii)}}]
$G^{\tau} P$ is closed in $G$.  
\end{enumerate}
If one of (hence all of ) the above three conditions is fulfilled then
${\mathfrak{p}}^{\tau}:={\mathfrak{p}} \cap {\mathfrak{g}}^{\tau}$
is a parabolic subalgebra of $\mathfrak{g}^\tau$,
and any irreducible subquotient occurring in the restriction
$X|_{\mathfrak {g}^{\tau}}$
 belongs to the category ${\mathcal{O}}^{\mathfrak {p}^{\tau}}$.
\end{theorem}
In Proposition \ref{prop:HPG}, 
 the geometric condition (iii)
 in Theorem \ref{thm:decoc} will be reformalised as an algebraic condition.

\medskip
\noindent
{\bf{Strategy of Proof of Theorem \ref{thm:decoc}:}}\enspace
We have already seen the equivalence (i) $\Longleftrightarrow$ (ii) in Lemma \ref{lem:YX}
 and the implication (iii) $\Longrightarrow$ (ii) in Proposition \ref{prop:closed}
 in a more general setting, i.e.\ without assuming that
$(\mathfrak{g},\mathfrak{g}')$ is a symmetric pair.  
The non-trivial part is the implication (ii) $\Longrightarrow$ (iii), 
 which will be proved in 
 Subsection \ref{subsec:pfdeco}
 after we establish some structural results on closed $G^{\tau}$-orbit
 in $G/P$
 (Subsection \ref{subsec:HPG}).  

\medskip

We end this subsection with two very special cases of
 Theorem \ref{thm:decoc}, namely,
for $\mathfrak{p} = \mathfrak{b}$ (Borel) and for the pair
$(\mathfrak{g} \oplus \mathfrak{g}, \operatorname{diag}\mathfrak{g})$:
\begin{cor}\label{cor:bdeco}
Let $\mathcal{O}$ be the BGG category associated to a Borel subalgebra
$\mathfrak{b}$, 
and $\tau$ an involutive automorphism of\/ $\mathfrak{g}$.
Then the following four conditions on
$(\tau,\mathfrak{b})$ are equivalent:
\begin{enumerate}[\upshape (i)]
\item  %(i)
Any simple $\mathfrak{g}$-module in $\mathcal{O}$ contains at least
one simple $\mathfrak{g}^\tau$-module when restricted to
$\mathfrak{g}^\tau$. 
\item  %(ii)
Any simple $\mathfrak{g}$-module in $\mathcal{O}$ is discretely
decomposable as a $\mathfrak{g}^\tau$-module in the sense of
Definition \ref{def:discdeco}.
\item  %(iii)
$G^\tau B$ is closed in $G$
\item  %(iv)
$\tau\mathfrak{b} = \mathfrak{b}$.
\end{enumerate}
\end{cor}

\begin{proof}
We shall see in Lemma \ref{lem:BGtau} that $G^\tau B$ is closed in $G$
if and only if $\tau\mathfrak{b} = \mathfrak{b}$.
Hence, Corollary follows from Theorem \ref{thm:decoc}.
\end{proof}

\begin{example}\label{ex:gl422}
Let $(\mathfrak{g},\mathfrak{g}')=(\mathfrak{sl}_4(\mathbb{C}),
\mathfrak{s}(\mathfrak{gl}_2(\mathbb{C})\oplus\mathfrak{gl}_2(\mathbb{C})))$.
Then there are $21$ orbits of the subgroup 
$S(GL(2,\mathbb{C})\times GL(2,\mathbb{C}))$ on the full flag variety
$\mathfrak{B}$ of 
$SL(4,\mathbb{C})$,
and $6$ $(=4!/2!2!)$ closed orbits among them.
In the following diagram,
 vertices stand for $21$ orbits on $\mathfrak{B}$,
 and edges generate their closure relations.
\begin{figure}[H]
\begin{center}
\psset{unit=2pt,linewidth=.3pt}
\pspicture(0,-45)(120,0)

	\rput(0,0){$\bullet$}
	\rput(20,0){$\bullet$}
	\rput(40,0){$\bullet$}
	\rput(60,0){$\bullet$}
	\rput(80,0){$\bullet$}
	\rput(100,0){$\bullet$}
	
	\pscircle[linecolor=red,linewidth=1pt](0,0){5}
	\pscircle[linecolor=red,linewidth=1pt](20,0){5}
	\pscircle[linecolor=red,linewidth=1pt](40,0){5}
	\pscircle[linecolor=red,linewidth=1pt](60,0){5}
	\pscircle[linecolor=red,linewidth=1pt](80,0){5}
	\pscircle[linecolor=red,linewidth=1pt](100,0){5}

	\psline{-}(0,-3)(0,-12)
	\psline{-}(20,-3)(0,-12)
	\psline{-}(20,-3)(20,-12)
	\psline{-}(20,-3)(38,-12)	%X2 - Y3
	\psline{-}(40,-3)(22,-12)	%X3 - Y2
	\psline{-}(40,-3)(60,-12)
	\psline{-}(60,-3)(42,-12)	%X4 - Y3
	\psline{-}(60,-3)(80,-12)
	\psline{-}(80,-3)(62,-12)	%X5 - Y4
	\psline{-}(80,-3)(80,-12)
	\psline{-}(80,-3)(100,-12)
	\psline{-}(100,-3)(100,-12)

	\rput(0,-15){$\bullet$}
	\rput(20,-15){$\bullet$}
	\rput(40,-15){$\bullet$}
	\rput(60,-15){$\bullet$}
	\rput(80,-15){$\bullet$}
	\rput(100,-15){$\bullet$}

	\psline{-}(0,-18)(8,-27)	%Y1 - Z1
	\psline{-}(0,-18)(30,-27)
	\psline{-}(18,-18)(10,-27)	%Y2 - Z1
	\psline{-}(20,-18)(48,-27)	%Y2 - W
	\psline{-}(38,-18)(32,-27)	%Y3 - Z2
	\psline{-}(42,-18)(49,-27)	%Y3 - W
	\psline{-}(58,-18)(51,-27)	%Y4 - W
	\psline{-}(60,-18)(67,-27)	%Y4 - Z3
	\psline{-}(80,-18)(52,-27)	%Y4 - W
	\psline{-}(80,-18)(90,-27)
	\psline{-}(100,-18)(70,-27)
	\psline{-}(100,-18)(92,-27)	%Y6 - Z4

	\rput(10,-30){$\bullet$}
	\rput(30,-30){$\bullet$}
	\rput(50,-30){$\bullet$}
	\rput(70,-30){$\bullet$}
	\rput(90,-30){$\bullet$}

	\psline{-}(12,-33)(23,-38)	%Z1 - a
	\psline{-}(29,-33)(26,-38)	%Z2 - a
	\psline{-}(50,-33)(50,-38)
	\psline{-}(70,-33)(74,-38)	%Z3 - c
	\psline{-}(88,-33)(77,-38)	%Z4 - c

	\rput(25,-40){$\bullet$}
	\rput(50,-40){$\bullet$}
	\rput(75,-40){$\bullet$}

	\psline{-}(27,-41)(50,-47)	%a - U
	\psline{-}(50,-42)(50,-47)
	\psline{-}(73,-41)(50,-47)	%c - U

	\rput(50,-50){$\bullet$}
	
	\rput(120,0){\textcolor{red}{closed}}
	\rput(120,-50){\textcolor{red}{open}}
	
\endpspicture
\end{center}
\caption{}
\label{fig:4.1}
\end{figure}
\noindent
Correspondingly, for a fixed Borel subalgebra $\mathfrak{b}$ of
$\mathfrak{g}$ there are\/ $6$ injective homomorphisms
$\iota_j: \mathfrak{g}'\to\mathfrak{g}$ $(1\le j\le 6)$
such that any simple $\mathfrak{g}$-module in $\mathcal{O}$ is discretely decomposable
when restricted to $\iota_j(\mathfrak{g}')$ and that
$\iota_j(\mathfrak{g}')$ is not conjugate to each other by an element
of the Borel subgroup.
\end{example}

\begin{cor}\label{cor:tensor}
Let\/ $\mathfrak{p}_1$, $\mathfrak{p}_2$ be two parabolic subalgebras of
a complex semisimple Lie algebra $\mathfrak{g}$.
Then the following three conditions on $(\mathfrak{p}_1,\mathfrak{p}_2)$
are equivalent:
\begin{enumerate}[\upshape (i)]
\item % (i)
For any simple $\mathfrak{g}$-module $X_1$ in
$\mathcal{O}^{\mathfrak{p}_1}$ and $X_2$ in
$\mathcal{O}^{\mathfrak{p}_2}$,
the tensor product representation
$X_1 \otimes X_2$ contains at least one simple $\mathfrak{g}$-module.
\item % (ii)
For any simple $\mathfrak{g}$-module $X_1$ in
$\mathcal{O}^{\mathfrak{p}_1}$ and $X_2$ in
$\mathcal{O}^{\mathfrak{p}_2}$,
the tensor product representation
$X_1 \otimes X_2$ is discretely decomposable as a
$\mathfrak{g}$-module.
\item % (iii)
$\mathfrak{p}_1 \cap \mathfrak{p}_2$ is a parabolic subalgebra.
\end{enumerate}
\end{cor}

\begin{proof}
Let $P_1$ and $P_2$ be the parabolic subgroups of
$G = \operatorname{Int}(\mathfrak{g})$
with Lie algebras $\mathfrak{p}_1$ and $\mathfrak{p}_2$, respectively.
Then the diagonal $G$-orbit on
$(G \times G) / (P_1 \times P_2)$ through the origin is given as
$G / (P_1 \cap P_2)$,
which is closed if and only if
$\mathfrak{p}_1 \cap \mathfrak{p}_2$ is a parabolic algebra of
$\mathfrak{g}$.
Hence, Corollary is deduced from Theorem \ref{thm:decoc}.
\end{proof}

\subsection{Criterion for closed $G^{\tau}$-orbit on $G/P$} 
\label{subsec:HPG}

As a preparation for the proof of Theorem \ref{thm:decoc},
we establish some structural results for closed $G^\tau$-orbits on the
generalized flag variety $G/P$ in this subsection.
We note that the closedness condition for $G^\tau$-orbits on
$G/P$ is much more
complicated than that for the full flag variety $G/B$
(cf.\ Lemma \ref{lem:BGtau} below).
The author is grateful to T. Matsuki for helpful discussions,
in particular,
for the proof of Proposition \ref{prop:HPG}.

Let $\mathfrak{g}$ be a complex semisimple Lie algebra,
$G = \operatorname{Int}(\mathfrak{g})$,
and $\tau$ an involutive automorphism of $\mathfrak{g}$
 as before.
We begin with:
\begin{lemma}
\label{lem:jtau}
\
\begin{enumerate}
\item[{\rm{1)}}]
Let $\theta$ be a Cartan involution of ${\mathfrak{g}}$
 commuting with $\tau$.  
For any parabolic subalgebra ${\mathfrak{p}}$, 
 there exist $h \in G^{\tau}$ and a Cartan subalgebra ${\mathfrak{j}}$
 such that $\tau {\mathfrak{j}} = \theta{\mathfrak{j}}={\mathfrak{j}}$
 and ${\mathfrak{j}} \subset \operatorname{Ad} (h) {\mathfrak{p}}$.  
In particular,
any parabolic subalgebra contains
a $\tau$-stable Cartan subalgebra.  
\item[{\rm{2)}}]
A parabolic subalgebra is $\tau$-stable 
 if and only if it is ${\mathfrak{g}}^{\tau}$-compatible
 (see Definition \ref{def:cp}).  
\end{enumerate}
\end{lemma}

\begin{proof}
{1)}\enspace
This assertion holds 
 for any Borel subalgebra of ${\mathfrak{g}}$ (\cite[Theorem 1]{M79}).  
Hence, 
 it holds also for any parabolic subalgebra.  
\newline
{2)}\enspace
Suppose ${\mathfrak{p}}$ is a $\tau$-stable parabolic subalgebra.  
Take a $\tau$-stable Cartan subalgebra ${\mathfrak{j}}$ 
 contained in ${\mathfrak{p}}$.  
Then there exists $H \in {\mathfrak{j}}$ such that
\[
  {\mathfrak{p}} = 
\bigoplus_{\substack{\alpha\in\Delta(\mathfrak{g},\mathfrak{j})\\
           \alpha(H)\ge 0}} 
{\mathfrak{g}}_{\alpha}.  
\]
Since $\tau{\mathfrak{p}} ={\mathfrak{p}}$,
$\alpha (H) \ge 0$
 if and only if $\alpha (\tau H) \ge 0$, 
 which is then equivalent
 to $\alpha(H + \tau H) \ge 0$.  
Therefore, the parabolic subalgebra
${\mathfrak{p}}$ equals ${\mathfrak{p}}(H+\tau H)$ 
with the notation \eqref{eqn:pH}, 
 and thus it is ${\mathfrak{g}}^{\tau}$-compatible.  
Conversely, 
 any ${\mathfrak{g}}^{\tau}$-compatible parabolic subalgebra 
 is obviously $\tau$-stable.  
\end{proof}

We then deduce a simple characterization of closed $G^{\tau}$-orbits 
 on the full flag variety $G/B$ from \cite[Proposition 2]{M79}
 combined with Lemma \ref{lem:jtau} 2):
\begin{lemma}
\label{lem:BGtau}
The following three conditions
 on $\tau$ and a Borel subalgebra ${\mathfrak{b}}$
 are equivalent:
\begin{enumerate}
\item[{\rm{(i)}}]
$G^{\tau}B$ is closed in $G$. 
\item[{\rm{(ii)}}]
$\tau{\mathfrak{b}} = {\mathfrak{b}}$.
\item[{\rm{(iii)}}]
${\mathfrak{b}}$ is ${\mathfrak{g}}^{\tau}$-compatible.  
\end{enumerate}
\end{lemma}
Unfortunately,
such a simple statement does not hold 
for a general parabolic subalgebra $\mathfrak{p}$.
In fact,
 the condition $\tau{\mathfrak{p}}={\mathfrak{p}}$
 is stronger than the closedness of $G^\tau P$ (see Example \ref{ex:cp}).  
In order to give the right characterization 
 for the closedness of $G^{\tau}P$, 
 we let $\operatorname{pr}_{\tau}:{\mathfrak{g}}\to {\mathfrak{g}}^{\tau}$
 be the projection
 defined by 
\begin{equation}
\label{eqn:prtau}
\operatorname{pr}_{\tau}(Z):=
  \frac 1 2 (Z+ \tau Z).  
\end{equation}

For a subspace $V$ in ${\mathfrak{g}}$, 
 we define the $\pm1$ eigenspaces of $\tau$ by
\begin{equation}
\label{eqn:Vtau}
V^{\pm\tau}:= \{v \in V : \tau v = \pm v\}.
\end{equation}  
Note that $\operatorname{pr}_{\tau}(V)=V^{\tau}$
 if $V$ is $\tau$-stable.  
\begin{proposition}
\label{prop:HPG}
Suppose ${\mathfrak {p}}$ is a parabolic subalgebra
 with nilradical ${\mathfrak{u}}$,
 and $\tau$ is an involutive automorphism of ${\mathfrak {g}}$.  
Then,
 the following three conditions
 on the triple $({\mathfrak{g}}, {\mathfrak{g}}^\tau, {\mathfrak{p}})$
 are equivalent:
\begin{enumerate}
\item[{\rm{(i)}}]$G^{\tau} P$ is closed in $G$.  
\item[{\rm{(ii)}}] 
 $\operatorname{pr}_{\tau}({\mathfrak {u}})$
 is a nilpotent Lie algebra.  
\item[{\rm{(iii)}}] 
 $\operatorname{pr}_{\tau}({\mathfrak {u}})$
 consists of nilpotent elements.    
\end{enumerate}
\end{proposition}

We note that the parabolic subalgebra $\mathfrak{p}$ may not be
$\tau$-stable in Proposition \ref{prop:HPG}.
The idea of the following proof goes back to \cite{M82},
which is to use a $\tau$-stable Borel subalgebra contained in
$\mathfrak{p}$ when $\mathfrak{p}$ itself is not $\tau$-stable.

\begin{proof}
We take a Borel subalgebra ${\mathfrak {b}} \subset {\mathfrak {p}}$
 such that $G^{\tau} B$ is relatively closed 
in $G^{\tau}P$.  
This is possible because $G^{\tau}\backslash G/B$
 is a finite set.

(i) $\Longrightarrow$ (ii) \enspace
Suppose $G^{\tau}P$ is closed in $G$.  
Then $G^{\tau} B$ is also closed in $G$.  
Owing to Lemma \ref{lem:BGtau},
 ${\mathfrak{b}}$ is $\tau$-stable,
and therefore, so is the nilradical 
 ${\mathfrak {n}}$ of ${\mathfrak {b}}$.  
Thus, 
 $\operatorname{pr}_{\tau}({\mathfrak{n}})={\mathfrak{n}}^{\tau}$.  
Since ${\mathfrak {u}} \subset {\mathfrak {n}}$, 
 we get 
$
  \operatorname{pr}_{\tau}({\mathfrak {u}}) 
  \subset 
  \operatorname{pr}_{\tau}({\mathfrak {n}})={\mathfrak {n}}^{\tau}$.  

For $X, Y \in {\mathfrak {g}}$, 
 a simple computation shows
\[
2 [\operatorname{pr}_{\tau}(X),
\operatorname{pr}_{\tau}(Y)]
=\operatorname{pr}_{\tau}([X,Y])
+
\operatorname{pr}_{\tau}([X, \tau Y]).  
\]
If $X, Y \in {\mathfrak {u}}$, 
then $[X,Y] \in {\mathfrak {u}}$
 and $[X, \tau Y] \in [{\mathfrak {u}}, {\mathfrak {n}}] \subset {\mathfrak {u}}$.  
Hence $\operatorname{pr}_{\tau}({\mathfrak {u}})$ is a Lie subalgebra.  
Since $\operatorname{pr}_{\tau}({\mathfrak {u}})$ is contained in ${\mathfrak {n}}^{\tau}$, 
 we conclude that $\operatorname{pr}_{\tau}({\mathfrak {u}})$
 is a nilpotent Lie algebra.    
Thus, 
 (i) $\Longrightarrow$ (ii) is proved.

(ii) $\Longrightarrow$ (iii).  Obvious.   

(iii) $\Longrightarrow$ (i). 
Since the conditions (i) and (iii) 
 remain the same 
 if we replace ${\mathfrak {p}}$
 by $\operatorname{Ad}(h) ({\mathfrak {p}})$
 for some $h \in G^{\tau}$, 
 we may and do assume
 that ${\mathfrak{p}}$ contains a Cartan subalgebra ${\mathfrak{j}}$
 such that $\tau {\mathfrak{j}}= \theta {\mathfrak{j}}= {\mathfrak{j}}$
 by Lemma \ref{lem:jtau}.  
Then $\theta \alpha = - \alpha$
 for any $\alpha\in \Delta({\mathfrak{g}}, {\mathfrak{j}})$.

Suppose $G^{\tau} P$ is not closed in $G$.
By the Matsuki duality \cite{M79},
 we see that $G^{\tau\theta} P$
 is not open in $G$.  
Therefore,
 there exists $\alpha \in \Delta({\mathfrak {u}}, {\mathfrak{j}})$
 such that ${\mathfrak {g}}_{-\alpha } \not \subset {\mathfrak {g}}^{\tau \theta} + \mathfrak {p}$.  
Take a non-zero $X_{-\alpha} \in {\mathfrak {g}}_{-\alpha}$.  
In view that 
\[
  X_{-\alpha}=(X_{-\alpha} + \tau \theta X_{- \alpha})- \tau \theta X_{-\alpha}
              \in {\mathfrak {g}}^{\tau\theta} + {\mathfrak {g}}_{\tau \alpha}, 
\]
we see 
${\mathfrak {g}}_{\tau \alpha} \not \subset {\mathfrak {p}}$
 because otherwise $X_{-\alpha}$ would be contained in ${\mathfrak {g}}^{\tau \theta} + {\mathfrak {p}}$.  
Hence,
$\mathfrak{g}_{-\tau\alpha} \subset \mathfrak{u}$
and 
$\tau \alpha \ne \alpha$.

Take a non-zero $X_{\alpha } \in {\mathfrak {g}}_{\alpha}$
 and we set $X:= X_{\alpha} + \tau \theta X_{\alpha} \in {\mathfrak {g}}_{\alpha} 
+{\mathfrak {g}}_{-\tau \alpha} \subset {\mathfrak {u}}$.  
\begin{description}
\item[\mdseries Case 1.]
Suppose $X \ne 0$. 
Let $Y:= \operatorname{pr}_{\tau} (X)$.  
Clearly, 
 $\theta Y = Y$.  
Moreover,
 $Y \ne 0$
 because $\tau \alpha \ne \alpha$.  
This means that $\operatorname{pr}_{\tau}({\mathfrak {u}})$
 contains a non-zero semisimple element.  
\item[\mdseries Case 2.]
Suppose $X=0$. 
Let $Y:= X_{\alpha}+ \tau X_{\alpha}= X_{\alpha} - \theta X_{\alpha}$.  
Then $Y \ne 0$
 and $\theta Y=-Y$.  
Again, this means that 
 $\operatorname{pr}_{\tau}({\mathfrak{u}})$
 contains a non-zero semisimple element.  
\end{description}

Thus we have proved the contraposition, 
 \lq\lq{not (i) $\Longrightarrow$ not (iii)}\rq\rq.  
Hence the proof of Proposition has been completed.
\end{proof}

The nilradical of the Lie algebra ${\mathfrak{p}}^{\tau}$
 is given explicitly as follows:

\begin{proposition}
\label{prop:ptau}
Under the equivalent conditions
 (i)-(iii) in Proposition \ref{prop:HPG},
 ${\mathfrak {p}}^{\tau}$
 is a parabolic subalgebra of ${\mathfrak {g}}^{\tau}$
 having the following Levi decomposition:
\[
   {\mathfrak {p}}^{\tau}
  ={\mathfrak {l}}^{\tau}
  +\operatorname{pr}_{\tau}({\mathfrak {u}}).
\]
\end{proposition}

\begin{proof}
[Proof of Proposition \ref{prop:ptau}]
We take a Borel subalgebra ${\mathfrak{b}}\subset{\mathfrak{p}}$
 such that $G^{\tau}B$
 is closed,
 and a $\tau$-stable Cartan subalgebra ${\mathfrak {j}}$
 contained in ${\mathfrak {b}}$
 as in the proof of Proposition \ref{prop:HPG}.

Given a ${\mathfrak {j}}$-stable subspace
 $V= \bigoplus_{\alpha \in \Delta(V)}  {\mathfrak {g}}_{\alpha}$
 in ${\mathfrak {g}}$, 
we denote by
$
   \Delta(V)
$
the multiset of\/ $\mathfrak{j}$-weights.
(Here we note that the multiplicity of the zero weight in $V$ may be
larger than one.) 
We divide
$\Delta(V)$
into the disjoint union
\[
    \Delta(V)=\Delta(V)_{\rm{I}} \amalg \Delta(V)_{\rm{II}}
  \amalg\Delta(V)_{\rm{III}}, 
\]
subject to the condition
{(I)}\enspace
$\tau \alpha = \alpha$
 and $\tau|_{\mathfrak {g}_{\alpha}}=\operatorname{id}$, 
{(II)}\enspace
$\tau \alpha = \alpha$
 and $\tau|_{{\mathfrak {g}}_{\alpha}}=-\operatorname{id}$, 
 and 
{(III)}\enspace
$\tau \alpha \ne \alpha$.
Accordingly, we have a direct sum as vector spaces:
\begin{align*}
    V^{\tau} =& \bigoplus_{\alpha \in \Delta(V)_{\operatorname{I}}}
                  {\mathfrak {g}}_{\alpha}
                \oplus
                 \bigoplus_{\alpha, \tau\alpha \in \Delta(V)_{\operatorname{III}}}
                 ({\mathfrak {g}}_{\alpha} + {\mathfrak {g}}_{\tau \alpha})^{\tau},
\\
   \operatorname{pr}_{\tau}(V)
             =& \bigoplus_{\alpha \in \Delta(V)_{\operatorname{I}}}
                  {\mathfrak {g}}_{\alpha}
                \oplus
                 \bigoplus_{\alpha \in \Delta(V)_{\operatorname{III}}}
                 ({\mathfrak {g}}_{\alpha} + {\mathfrak {g}}_{\tau \alpha})^{\tau}.  
\end{align*}
In particular, 
 we get 
\begin{align*}
{\mathfrak{p}}^{\tau}
=&\bigoplus_{\alpha \in \Delta({\mathfrak{p}})_{\operatorname{I}}} {\mathfrak{g}}_{\alpha}
   \oplus
   \bigoplus_{\alpha, \tau\alpha \in \Delta({\mathfrak{p}})_{\operatorname{III}}}
  ({\mathfrak{g}}_{\alpha} + {\mathfrak{g}}_{\tau \alpha})^{\tau}
\\
=&\bigoplus_{\alpha \in \Delta({\mathfrak{l}})_{\operatorname{I}}} {\mathfrak{g}}_{\alpha}
   \oplus
   \bigoplus_{\alpha, \tau\alpha \in \Delta({\mathfrak{l}})_{\operatorname{III}}}
   ({\mathfrak{g}}_{\alpha} + {\mathfrak{g}}_{\tau \alpha})^{\tau}
\oplus\bigoplus_{\alpha \in \Delta({\mathfrak{u}})_{\operatorname{I}}} {\mathfrak{g}}_{\alpha}
  \oplus
   \bigoplus_{\alpha, \tau\alpha \in \Delta({\mathfrak{u}})_{\operatorname{III}}}
   ({\mathfrak{g}}_{\alpha} + {\mathfrak{g}}_{\tau \alpha})^{\tau}
\\
=&{\mathfrak {l}}^{\tau} \oplus \operatorname{pr}_{\tau}({\mathfrak {u}}).  
\end{align*}
Here we have used $\tau {\mathfrak {u}} \subset {\mathfrak {p}}$
in the second equality.  
Thus Proposition \ref{prop:ptau} is proved.  
\end{proof}

\subsection{Application of associated varieties to restrictions}
\label{subsec:AV}

In this subsection,
we apply 
 associated varieties of $\mathfrak{g}$-models
to the study of branching problems.

Suppose $X$ is a finitely generated ${\mathfrak {g}}$-module.  
We take a finite dimensional subspace $X_0$
 which generates $X$ as a ${\mathfrak {g}}$-module.  
Let 
$
  U({\mathfrak{g}})= \cup_{k \ge 0}U_k({\mathfrak{g}})
$
 be a natural filtration 
 of the enveloping algebra of ${\mathfrak{g}}$. 
Then, $X_k:=U_k({\mathfrak {g}}) X_0$
 ($k \in {\mathbb{N}}$)
 gives a filtration $\{X_k\}_k$
satisfying
\[
     X= \bigcup_{k=0}^{\infty}X_k,
    \qquad
    U_i({\mathfrak {g}}) X_j = X_{i+j}\,\,
    (i,j \ge 0). 
\] 
Then,
 $\operatorname{gr}X:=\bigoplus_{k=0}^{\infty} X_k/X_{k-1}$ is a
 finitely generated module of the commutative algebra 
$
    \operatorname{gr}U({\mathfrak {g}})\simeq S({\mathfrak {g}})
$.  
The {\it{associated variety}}
 of the ${\mathfrak {g}}$-module $X$ is a closed subset ${\mathcal{V}}_{\mathfrak{g}}(X)$
 of ${\mathfrak{g}}^{\ast}$
 defined by 
\[
  {\mathcal{V}}_{\mathfrak {g}}(X)
  :=\operatorname{Supp}_{S({\mathfrak {g}})}(\operatorname{gr}X).  
\] 
Then ${\mathcal{V}}_{\mathfrak {g}}(X)$ is independent
 of the choice of the generating subspace $X_0$.  
We recall the following basic properties:

\begin{lemma}
[{\cite[Chapter 17]{Ja}}]
\label{lem:AV}
{\rm{1)}}\enspace
If $0\longrightarrow X_1 \longrightarrow X \longrightarrow X_2\longrightarrow 0$
 is an exact sequence of ${\mathfrak{g}}$-modules,
 we have ${\mathcal{V}}_{\mathfrak{g}}(X)={\mathcal{V}}_{\mathfrak{g}}(X_1) \cup 
 {\mathcal{V}}_{\mathfrak{g}}(X_2)$.  
\par\noindent
{\rm{2)}}\enspace
For any finite dimensional ${\mathfrak {p}}$-module $F$, 
$
     {\mathcal{V}}_{\mathfrak {g}}(U({\mathfrak {g}}) \otimes_{U({\mathfrak {p}})}F)
     = 
     {\mathfrak {p}}^{\perp}
$.
\end{lemma}
Let ${\mathfrak {g}}'$ be a reductive subalgebra in ${\mathfrak {g}}$,
and
  $\operatorname{pr}_{{\mathfrak {g}} \to {\mathfrak {g}}'}: {\mathfrak {g}}^* \to {{\mathfrak {g}}'}^*$
 the restriction map.     
We set 
$
   {\mathfrak{p}}':= {\mathfrak{g}}' \cap {\mathfrak{p}}
$
  and 
$
   {\mathfrak{p}}'^{\perp}
   :=
   \{\lambda \in ({\mathfrak{g}}')^{\ast}:
     \lambda|_{{\mathfrak{p}}'}\equiv 0
   \}.
$
\begin{lemma}
\label{lem:AVd}
Let $X$ be a simple 
${\mathfrak{g}}$-module lying in ${\mathcal{O}}^{\mathfrak{p}}$,
and $\mathfrak{g}'$ a reductive subalgebra in $\mathfrak{g}$.

\noindent
{\rm{1)}}\enspace
If $Y$ is a simple ${\mathfrak {g}}'$-module
 such that $\operatorname{Hom}_{\mathfrak {g}'}(Y,X|_{\mathfrak {g}'}) \ne \{0\}$
 then 
\begin{equation}
\label{eqn:VY}
  \operatorname{pr}_{{\mathfrak {g}}\rightarrow {\mathfrak {g}}'}
     ({\mathcal{V}}_{\mathfrak {g}}(X))
    \subset
    {\mathcal{V}}_{\mathfrak {g}'}(Y) \subset ({\mathfrak{p}}')^{\perp}.  
\end{equation}

\par\noindent
{\rm{2)}}\enspace 
If $Y_i$ are simple ${\mathfrak{g}}'$-modules
 such that $\operatorname{Hom}_{\mathfrak{g}'}(Y_i,X|_{\mathfrak{g}'})\ne \{0\}$
 ($i=1,2$), 
then
$
  {\mathcal{V}}_{\mathfrak{g}'}(Y_1)= {\mathcal{V}}_{\mathfrak{g}'}(Y_2)
$.
\end{lemma} 
\begin{proof}
1)\enspace 
Since ${\mathcal{O}}^{\mathfrak {p}}$ is closed under tensor products
 with finite dimensional representations,
 the proof for the first inclusion in \eqref{eqn:VY}
 parallels to the proof of \cite[Theorem 3.1]{Ko3} 
 by using the double filtration of $X$.

For the second inclusion in \eqref{eqn:VY}, 
 we use the notation of the proof of Proposition \ref{prop:closed}
 and let $Y$ be the image of \eqref{eqn:YXlmd}.
Then it follows from Lemma \ref{lem:AV}
 that
${\mathcal{V}}_{\mathfrak {g}'}(Y)
    \subset
    {\mathcal{V}}_{\mathfrak {g}'}
(U({\mathfrak {g}}') \otimes_{U({\mathfrak {p}}')}
     (F_{\lambda}|_{\mathfrak {p}'}))
={\mathfrak {p}}'^{\perp}.  
$

\par\noindent
{2)} \enspace
The proof is the same as that of \cite[Theorem 3.7]{Ko3}
 in the category of $({\mathfrak{g}},K)$-modules.  
\end{proof}

\begin{remark}\label{rem:AVd}
An analogous result to Lemma \ref{lem:AVd} 2) was shown in \cite{xEW}
 in the special case where $X$ is the oscillator representation
 of ${\mathfrak{g}}={\mathfrak{sp}}(n,{\mathbb{R}})$
 in the context of compact dual pair correspondence by case-by-case computations.  
In this case,
$\operatorname{pr}_{{\mathfrak {g}} \to {\mathfrak {g}}'}
 (\mathcal{V}_{\mathfrak{g}}(X))$
coincides with the associated variety
$\mathcal{V}_{\mathfrak{g}'}(Y)$.
It is plausible that
$\operatorname{pr}_{{\mathfrak {g}} \to {\mathfrak {g}}'}
 (\mathcal{V}_{\mathfrak{g}}(X))
 =  \mathcal{V}_{\mathfrak{g}'}(Y)$
in the generality of the setting in Lemma \ref{lem:AVd}.
We shall discuss this assertion in
 Theorem \ref{thm:AV} below for symmetric pairs
$(\mathfrak{g},\mathfrak{g}^\tau)$.
\end{remark}

\subsection{Proof of Theorem \ref{thm:decoc}}
\label{subsec:pfdeco}

The equivalence of Theorem \ref{thm:decoc} has been
 already proved  in Section \ref{subsec:2.4}
 except for the implication (ii) $\Rightarrow$ (iii).  
We are ready to complete the proof.   

\begin{proof}
[Proof of Theorem \ref{thm:decoc}, 
(ii) $\Rightarrow$ (iii)]
By the Killing form,
 we identify ${\mathfrak{g}}^*$ with ${\mathfrak{g}}$.  
Then the projection $\operatorname{pr}_{{\mathfrak {g}} \to {\mathfrak {g}}^{\tau}}:
 {\mathfrak{g}}^* \to ({\mathfrak{g}}^{\tau})^*$
 is given as the map 
$\operatorname{pr}_{\tau}: {\mathfrak {g}} \to {\mathfrak {g}}^{\tau}$
 (see \eqref{eqn:prtau}).  
Further,
 ${\mathfrak {p}}^{\perp} = \{ \lambda \in \mathfrak{g}^*:
   \lambda|_{\mathfrak{p}} \equiv 0\}$
 is isomorphic to the nilpotent radical $\mathfrak {u}$ of the parabolic subalgebra ${\mathfrak{p}}$.

We take a generalized Verma module 
$X := \ensub{\mathfrak{g}}{\lambda}{\mathfrak{p}}$
with generic parameter
$\lambda \in \Lambda^+(\mathfrak{l})$ (cf.\ \eqref{eqn:anti}).  
Then it follows from Lemma \ref{lem:AV}
 that ${\mathcal{V}}_{\mathfrak{g}}(X)={\mathfrak{u}}$.  
Therefore,
 if the restriction $X|_{\mathfrak{g}^{\tau}}$
 is discretely decomposable, 
 then $\operatorname{pr}_{\tau}({\mathfrak{u}})$
 consists of nilpotent elements by Lemma \ref{lem:AVd}.  
In turn,
 $G^{\tau} P$ is closed in $G$ owing to Proposition \ref{prop:HPG}.  
Thus, 
 the proof of Theorem \ref{thm:decoc} is completed.  
\end{proof}

\subsection{Associated varieties of irreducible summands}
\label{subsec:AVthm}

We retain the previous notation:
$\mathfrak{p}$ is a parabolic subalgebra of a complex semisimple Lie
algebra $\mathfrak{g}$,
and $\tau$ an involutive automorphism of $\mathfrak{g}$.
In this subsection,
we give an explicit formula
 for the associated variety $\mathcal{V}_{\mathfrak{g}^\tau}(Y)$
and the Gelfand--Kirillov dimension $\operatorname{DIM}(Y)$ 
 of irreducible summands $Y$.
\begin{theorem}
\label{thm:AV}
Suppose $({\mathfrak {g}}, {\mathfrak{g}}^\tau, {\mathfrak {p}})$ satisfies
 one of $($hence, all of$)$ the equivalent conditions
 in Theorem \ref{thm:decoc}.  
Let $X = \ensub{\mathfrak{g}}{\lambda}{\mathfrak{p}}$ 
be a simple generalized Verma module, 
 and $Y$ a simple ${\mathfrak{g}}'$-module
 such that 
$
  \operatorname{Hom}_{\mathfrak{g}'}
  (Y,X|_{\mathfrak{g}'}) \ne\{0\}.
$
Then,
\begin{equation*}
     {\mathcal{V}}_{\mathfrak {g}^{\tau}}(Y) = \operatorname{pr}_{\tau}({\mathfrak {u}})
\quad
\text{and}
\quad
   \operatorname{DIM}(Y) = \dim {\mathfrak {g}}^{\tau}/{\mathfrak {p}}^{\tau}.  
\end{equation*}
\end{theorem}

\begin{proof}
[Proof of Theorem \ref{thm:AV}]
The nilradical of the parabolic subalgebra
 ${\mathfrak{p}}^{\tau}$ is given by
 $\operatorname{pr}_{\tau} ({\mathfrak{u}})$
in Proposition \ref{prop:ptau}.  
Hence, via the isomorphism ${\mathfrak{g}}^* \simeq {\mathfrak{g}}$, 
the inclusive relation \eqref{eqn:VY}
 is written as 
\begin{equation}
\label{eqn:VY2}
   \operatorname{pr}_{\tau} ({\mathcal{V}}_{\mathfrak {g}}(X))
  \subset 
  {\mathcal{V}}_{\mathfrak {g}^{\tau}}(Y)
  \subset 
  \operatorname{pr}_{\tau} ({\mathfrak {u}}).
\end{equation}
Since ${\mathcal{V}}_{\mathfrak {g}}(X)={\mathfrak {u}}$,
 the three terms in \eqref{eqn:VY2} must be the same,
 and therefore
${\mathcal{V}}_{\mathfrak {g}^{\tau}}(Y) = \operatorname{pr}_{\tau}({\mathfrak{u}})$.

The Gelfand--Kirillov dimension ${\operatorname{DIM}}(Y)$
 is given by
the dimension of the associated variety
${\mathcal{V}}_{\mathfrak {g}^{\tau}}(Y)$, 
and thus we have 
${\operatorname{DIM}}(Y)=\operatorname{dim} \operatorname{pr}_{\tau}({\mathfrak{u}})$, 
 which equals 
$
\operatorname{dim}{\mathfrak{p}}^{\tau}
-
\operatorname{dim}{\mathfrak{l}}^{\tau}
=
\operatorname{dim}{\mathfrak{g}}^{\tau}
-
\operatorname{dim}{\mathfrak{p}}^{\tau}
$
 by Proposition \ref{prop:ptau}.  
\end{proof}

\begin{remark}
\label{rem:AV}

There are finitely many $G^{\tau}$-orbits on the generalized flag
variety $G/P$ by \cite{M79}.  
Among them, 
 suppose $G^{\tau}y_j P$
 ($j=1,2,\cdots, k$)
 are closed in $G$.  
Correspondingly
 we realize ${\mathfrak{g}}^{\tau}$
 as a subalgebra of ${\mathfrak{g}}$
 by 
\[
  \iota_j: {\mathfrak{g}}^{\tau} \hookrightarrow
{\mathfrak{g}},
 \quad
  Z \mapsto \operatorname{Ad}(y_j)^{-1}(Z).  
\]
Then $({\mathfrak{g}},\iota_j({\mathfrak{g}}^{\tau}))$
form symmetric pairs
 defined by the involutions
$
  \tau_j := \operatorname{Ad}(y_j^{-1}) \circ \tau \circ
\operatorname{Ad}(y_j)
\in \operatorname{Aut}({\mathfrak{g}}).  
$
Theorem \ref{thm:decoc} implies that
 the restrictions $X|_{\iota_j({\mathfrak{g}}^{\tau})}$
 are discretely decomposable for any $X \in {\mathcal{O}}^{\mathfrak{p}}$
and for any $j$ ($j=1,\dots,k$).  
Obviously, 
the Lie algebras $\iota_j(\mathfrak {g}^{\tau})$
 are isomorphic to each other,
 but $\dim ({\mathfrak {p}} \cap \iota_j ({\mathfrak {g}}^{\tau}))$
 may differ.
Accordingly,
 the Gelfand--Kirillov dimension
 of simple summands
 in the restrictions $X|_{\iota_j({\mathfrak{g}}^{\tau})}$
 depends on $j$.  
See Examples \ref{ex:Heisen1} and \ref{ex:GKCA} below.  
\end{remark}

\begin{example}
[$A_{p+q-1} \downarrow A_{p-1} \times A_{q-1}$]
\label{ex:Heisen1}
Let $p,q \ge 2$, ${\mathfrak{g}} = {\mathfrak{sl}}_{p+q}({\mathbb{C}})$, 
 ${\mathfrak{p}}$ its parabolic subalgebra
 whose nilradical is the Heisenberg Lie algebra of dimension $2(p+q)-3$,  
 and\/ ${\mathfrak{g}}'={\mathfrak{s}}({\mathfrak{gl}}_p({\mathbb{C}})\oplus {\mathfrak{gl}}_q({\mathbb{C}}))$.
Then,
 there are four injective homomorphisms
$
  \iota_j:{\mathfrak{g}}' \to {\mathfrak{g}}
$
 $(1 \le j \le 4)$
 such that each $\iota_j$ induces closed $G'$-orbits on $G/P$
and that $\iota_j(\mathfrak{g}')$ is not conjugate to each other by an
 element of $P$,
the parabolic subgroup with Lie algebra $\mathfrak{p}$.

The following diagram for $p=q=2$ shows how the $21$ orbits 
of $G'$ on the full flag variety $\mathfrak{B}\simeq G/B$ (see
 Figure
 \ref{fig:4.1}) are sent to $G'\backslash G/P$ under the quotient map
\[
G'\backslash G/B \to G'\backslash G/P
\]
induced from the inclusion $B\subset P$.
In particular, there are $4$ closed $G'$-orbits on $G/P$ among $10$
orbits.

\definecolor{darkgray}{gray}{0.5}

\psset{unit=1.6pt,linewidth=.6pt,hatchwidth=.4pt,hatchsep=2pt,hatchcolor=darkgray}

\begin{figure}[H]
\hbox{%
\begin{tabular}[m]{c}
\pspicture(-20,-60)(100,10)

	\psline[linecolor=red]{-}(0,0)(0,-15)
	\psline[linecolor=red]{-}(20,0)(0,-15)
	\psline[linecolor=blue]{-}(20,0)(20,-15)
	\psline[linecolor=blue]{-}(40,0)(20,-15)	%X3 - Y2
	\psline[linecolor=green]{-}(20,0)(40,-15)
	\psline[linecolor=green]{-}(40,0)(60,-15)
	\psline[linecolor=green]{-}(60,0)(40,-15)	%X4 - Y3
	\psline[linecolor=green]{-}(80,0)(60,-15)
	\psline[linecolor=blue]{-}(60,0)(80,-15)
	\psline[linecolor=blue]{-}(80,0)(80,-15)
	\psline[linecolor=red]{-}(80,0)(100,-15)
	\psline[linecolor=red]{-}(100,0)(100,-15)

	\psline[linecolor=blue]{-}(0,-15)(10,-30)	%Y1 - Z1
	\psline[linecolor=red]{-}(20,-15)(10,-30)	%Y2 - Z1
	\psline[linecolor=red]{-}(40,-15)(30,-30)
	\psline[linecolor=green]{-}(0,-15)(30,-30)	%Y2 - W
	\psline[linecolor=green]{-}(20,-15)(50,-30)	%Y2 - W
	\psline[linecolor=blue]{-}(40,-15)(50,-30)	%Y3 - W
	\psline[linecolor=blue]{-}(60,-15)(50,-30)	%Y4 - W
	\psline[linecolor=green]{-}(80,-15)(50,-30)	%Y4 - W
	\psline[linecolor=green]{-}(100,-15)(70,-30)	%Y2 - W
	\psline[linecolor=red]{-}(60,-15)(70,-30)
	\psline[linecolor=red]{-}(80,-15)(90,-30)
	\psline[linecolor=blue]{-}(100,-15)(90,-30)	%Y6 - Z4

	\psline[linecolor=green]{-}(10,-30)(25,-40)	%Z1 - a
	\psline[linecolor=blue]{-}(30,-30)(25,-40)	%Z2 - a
	\psline[linecolor=red]{-}(50,-30)(50,-40)
	\psline[linecolor=blue]{-}(70,-30)(75,-40)	%Z3 - c
	\psline[linecolor=green]{-}(90,-30)(75,-40)	%Z4 - c

	\psline[linecolor=red]{-}(25,-40)(50,-50)	%a - U
	\psline[linecolor=green]{-}(49.7,-40)(49.7,-50)
	\psline[linecolor=blue]{-}(50.3,-40)(50.3,-50)
	\psline[linecolor=red]{-}(75,-40)(50,-50)	%c - U

	\rput(0,0){$\bullet$}
	\rput(20,0){$\bullet$}
	\rput(40,0){$\bullet$}
	\rput(60,0){$\bullet$}
	\rput(80,0){$\bullet$}
	\rput(100,0){$\bullet$}
	
	\rput(0,-15){$\bullet$}
	\rput(20,-15){$\bullet$}
	\rput(40,-15){$\bullet$}
	\rput(60,-15){$\bullet$}
	\rput(80,-15){$\bullet$}
	\rput(100,-15){$\bullet$}

	\rput(10,-30){$\bullet$}
	\rput(30,-30){$\bullet$}
	\rput(50,-30){$\bullet$}
	\rput(70,-30){$\bullet$}
	\rput(90,-30){$\bullet$}

	\rput(25,-40){$\bullet$}
	\rput(50,-40){$\bullet$}
	\rput(75,-40){$\bullet$}

	\rput(50,-50){$\bullet$}
	
%	\rput(-15,0){closed}
%	\rput(-15,-50){open}
	
	\rput(50,-60){\( G' \backslash G / B \)}

	\pscircle[fillstyle=hlines,linecolor=darkgray](40,0){5}
	\pscircle[fillstyle=hlines,linecolor=darkgray](60,0){5}
	\pspolygon[linearc=3,fillstyle=hlines,linecolor=darkgray](-3,3)(-3,-20)(30,3)
	\pspolygon[linearc=3,fillstyle=hlines,linecolor=darkgray](103,3)(103,-20)(70,3)
	\rput{330}(15,-22.5){\psellipse[fillstyle=hlines,linecolor=darkgray](0,0)(5,12)}
	\rput{330}(35,-22.5){\psellipse[fillstyle=hlines,linecolor=darkgray](0,0)(5,12)}
	\rput{30}(65,-22.5){\psellipse[fillstyle=hlines,linecolor=darkgray](0,0)(5,12)}
	\rput{30}(85,-22.5){\psellipse[fillstyle=hlines,linecolor=darkgray](0,0)(5,12)}
	\psellipse[fillstyle=hlines,linecolor=darkgray](50,-35)(4,7)
	\psccurve[fillstyle=hlines,linecolor=darkgray](50,-43)(25,-37)(22,-40)(50,-53)(78,-40)(75,-37)
	
\endpspicture
\end{tabular}
\begin{tabular}[m]{c}
$\to$
\end{tabular}
\begin{tabular}[m]{c}
%\psset{unit=2pt,linewidth=.6pt,hatchwidth=.4pt,hatchsep=2pt,hatchcolor=darkgray}
\pspicture(0,-60)(120,10)

	\psline{-}(0,0)(10,-30)
	\psline{-}(40,0)(10,-30)
	\psline{-}(40,0)(70,-30)
	\psline{-}(60,0)(30,-30)
	\psline{-}(60,0)(90,-30)
	\psline{-}(100,0)(90,-30)

	\psline{-}(10,-30)(50,-40)
	\psline{-}(30,-30)(50,-40)
	\psline{-}(70,-30)(50,-40)
	\psline{-}(90,-30)(50,-40)

	\psline{-}(50,-40)(50,-50)

	\rput(0,0){$\bullet$}
%	\rput(20,0){$\bullet$}
	\rput(40,0){$\bullet$}
	\rput(60,0){$\bullet$}
%	\rput(80,0){$\bullet$}
	\rput(100,0){$\bullet$}
	
	\rput(0,5){$\iota_1$}
	\rput(40,5){$\iota_2$}
	\rput(60,5){$\iota_3$}
	\rput(100,5){$\iota_4$}
	
%	\rput(0,-15){$\bullet$}
%	\rput(20,-15){$\bullet$}
%	\rput(40,-15){$\bullet$}
%	\rput(60,-15){$\bullet$}
%	\rput(80,-15){$\bullet$}
%	\rput(100,-15){$\bullet$}

	\rput(10,-30){$\bullet$}
	\rput(30,-30){$\bullet$}
%	\rput(50,-30){$\bullet$}
	\rput(70,-30){$\bullet$}
	\rput(90,-30){$\bullet$}

%	\rput(25,-40){$\bullet$}
	\rput(50,-40){$\bullet$}
%	\rput(75,-40){$\bullet$}

	\rput(50,-50){$\bullet$}
	
%	\rput(115,0){closed}
%	\rput(115,-50){open}
	
	\rput(50,-60){\( G' \backslash G / P \)}
\endpspicture
\end{tabular}
}
\caption{}
\label{fig:4.2}
\end{figure}

For general $p,q\ge2$,
the number of closed $G'$-orbits on $G/P$ remains to be four.
Let $\tau_j \in \operatorname{Aut}({\mathfrak{g}})$
$(j=1,2,3,4)$
be defined as in Remark \ref{rem:AV}.  
It turns out that ${\mathfrak {p}}$ is $\tau_j$-compatible
 for all $j$.
Further, 
by applying Theorem \ref{thm:AV},
we see that
 the Gelfand--Kirillov dimension is given by 
\begin{equation*}
\operatorname{DIM}(Y)=
\begin{cases}
p+q-2 \qquad&(j=2,3)
\\
2p-3 &(j=1)
\\
2q-3&(j=4)  
\end{cases}
\end{equation*}
   for any simple ${\mathfrak{g}}'$-module $Y$
 and for any simple generalized Verma module $X =\ensub{\mathfrak{g}}{\lambda}{\mathfrak{p}}$
 such that
 $\operatorname{Hom}_{\mathfrak{g}'}(Y,X|_{\iota_j({\mathfrak{g}}')})\ne
 0$. 
By using Corollary \ref{cor:br},
we can find the branching laws of the
 restriction
$M_{\mathfrak{p}}^{\mathfrak{g}}(\lambda)|_{\iota_j(\mathfrak{g}')}$ 
for generic $\lambda$.
They are
 multiplicity-free for any $j=1,2,3,4$.
\end{example}

\begin{example}
[$C_n \downarrow A_n$]
\label{ex:GKCA}
Let ${\mathfrak {g}}$ be the complex
 symplectic Lie algebra $\mathfrak {sp}_n({\mathbb{C}})$ of rank $n$, 
 ${\mathfrak {p}}$ the Siegel parabolic subalgebra,
 and ${\mathfrak {g}}'={\mathfrak {gl}}_n({\mathbb{C}})$.  
Then there are $(n+1)$ injective homomorphisms
$\iota_j:{\mathfrak{g}}' \to {\mathfrak{g}}$
 ($0 \le j \le n$)
 such that 
$\iota_j(\mathfrak{g}')$ is not conjugate to each other under the
 Siegel parabolic subgroup
each $\iota_j$ induces closed
$GL_n(\mathbb{C})$-orbits on
$Sp(n,\mathbb{C})/P$ and that
\[
  \operatorname{DIM}(Y)=j(n-j)
\]
for any simple ${\mathfrak {g}}'$-module 
 if $Y$ occurs in the restriction $X|_{\iota_j({\mathfrak {g}}')}$
where $X$ is
 any simple generalized Verma module $\ensub{\mathfrak{g}}{\lambda}{\mathfrak{p}}$.   
\end{example}

\begin{proof}
[Sketch of the proof]
We take $\iota_j$ so that 
${\mathfrak {g}}_{\mathbb{R}} \simeq \mathfrak{sp}(n,\mathbb{R})$
and 
${\mathfrak {g}}_{\mathbb{R}} \cap \iota_j({\mathfrak{g}}') \simeq
{\mathfrak {u}}(j,n-j)$ 
where the specific real form $\mathfrak{g}_{\mathbb{R}}$
 will be explained in
Subsection \ref{subsec:PAbel}.
\end{proof}

\subsection{Finite multiplicity theorem}
\label{subsec:fmthm}

The multiplicities in branching laws
 behave much mildly in the BGG category ${\mathcal{O}}$
 than those 
 in the context of unitary representations
 (see Example \ref{ex:psmult} below).

Here is a finite multiplicity theorem in the category $\mathcal{O}$.

\begin{theorem}[finite multiplicity theorem]
\label{thm:fm}
Let $\tau$ be an involutive automorphism of a complex semisimple Lie
algebra $\mathfrak{g}$,
 and  $\mathfrak{b}$ a $\tau$-stable
 Borel subalgebra of $\mathfrak{g}$.
Then
\begin{equation}\label{eqn:MFtau}
\dim \operatorname{Hom}_{\mathfrak{g}^\tau}
(Y, X|_{\mathfrak{g}^\tau}) < \infty
\end{equation}
for any simple $\mathfrak{g}$-module $X$ 
in the category $\mathcal{O}\equiv \mathcal{O}^\mathfrak{b}$
and any simple $\mathfrak{g}^\tau$-module $Y$.
\end{theorem}

\begin{proof}
[Proof of Theorem \ref{thm:fm}]
Let $Y$ be any simple $\mathfrak{g}^\tau$-module.
We apply  
 Theorem \ref{thm:br} to the ${\mathfrak{g}}^{\tau}$-compatible Borel subalgebra ${\mathfrak{b}}$, 
 and conclude that
 the multiplicities of $Y$ occurring
 as subquotients of the restriction of any
 Verma module
$\ensub{\mathfrak{g}}{\lambda}{\mathfrak{b}}$ are finite.
Since any simple $\mathfrak{g}$-module $X \in \mathcal{O}$ is obtained
as the subquotient of some Verma module,
\eqref{eqn:MFtau} follows.
\end{proof}

\begin{remark}\label{rem:fm}
We recall that Theorem \ref{thm:br} counts the multiplicities in the
subquotients.
Therefore,
the multiplicities of $Y$ occurring in the restriction
$X|_{\mathfrak{g}^\tau}$
as \textit{subquotients} are also finite.
\end{remark}

Theorem \ref{thm:fm} should be compared with the fact
 that the multiplicities are often infinite
 in the branching laws of the restriction of an irreducible unitary representation
 with respect to a semisimple symmetric pair
(see \cite{Ko07}): 

\begin{example}
\label{ex:psmult}
There exists an irreducible unitary representation $\pi$ of $G=SO(5,{\mathbb{C}})$
 and two irreducible unitary representations $Y_1$ and $Y_2$ of the
 subgroup $G'=SO(3,2)$
 satisfying the following three conditions:
\begin{enumerate}
\item[{\rm{(1)}}]
$0 < \dim \operatorname{Hom}_{G'}
 (Y_1, \pi|_{G'}) < \infty$.
\item[{\rm{(2)}}]
$\dim \operatorname{Hom}_{G'}
 (Y_2, \pi|_{G'}) = \infty$.
\item[{\rm{(3)}}]
$\operatorname{DIM}(Y_1)=3$, 
$\operatorname{DIM}(Y_2)=4$.  
\end{enumerate}
Here, 
$\operatorname{Hom}_{G'}(\cdot, \cdot)$ denotes the space of
continuous $G'$-intertwining operators,
and
$\operatorname{DIM}(Y)$ stands for the Gelfand--Kirillov dimension of
the underlying $(\mathfrak{g}',K')$-module of the unitary
representation $Y$ of $G'$.
\end{example}

\section{Multiplicity-free branching laws}
\label{sec:mf}

In this section we prove two multiplicity-free theorems
 for the restriction of generalized Verma modules
 with respect to symmetric pairs $({\mathfrak{g}}, {\mathfrak{g}}')$:

$\bullet$\enspace
$\mathfrak{p}$ is special and $({\mathfrak{g}}, {\mathfrak{g}}')$ is general (Theorem \ref{thm:PAbel}), 

\medskip
$\bullet$\enspace
$\mathfrak{p}$ is general and $({\mathfrak{g}}, {\mathfrak{g}}')$ is special (Theorem \ref{thm:mfall}).  
\vskip 1pc
\par\noindent
Correspondingly, 
 explicit branching laws are also derived 
 (Theorems \ref{thm:mfbr1}, \ref{thm:AA}, \ref{thm:BD}, and \ref{thm:DB}).  

\subsection{Parabolic subalgebra with abelian nilradical}
\label{subsec:PAbel}

We begin with multiplicity-free branching laws of the restriction
$\ensub{\mathfrak{g}}{\lambda}{\mathfrak{p}}|_{\mathfrak{g}^\tau}$
with respect to symmetric pairs
$(\mathfrak{g},\mathfrak{g}^\tau)$
in the case 
where $\mathfrak{p}$ is a certain maximal parabolic subalgebra.

An abstract feature of the results here boils down to the following:
\begin{theorem}
\label{thm:PAbel}
Suppose ${\mathfrak {p}} = \mathfrak{l} + \mathfrak{u}_+$ 
is a parabolic subalgebra
such that the nilradical $\mathfrak{u}_+$ is abelian.
Then for any involutive automorphism $\tau$ such that $\tau {\mathfrak {p}}={\mathfrak {p}}$, 
 the generalized Verma module $\ensub{\mathfrak{g}}{\lambda}{\mathfrak{p}}$ of scalar type
 is decomposed into a multiplicity-free direct sum
 of simple ${\mathfrak {g}}^{\tau}$-modules
if $\lambda \in \Lambda^+(\mathfrak{l})$ 
 is sufficiently negative,
 i.e. $\langle \lambda, \alpha \rangle \ll 0$
 for all $\alpha \in \Delta({\mathfrak{u}}_+)$.  
\end{theorem}

Theorem \ref{thm:PAbel} is deduced from an explicit formula
 of the irreducible decomposition.
To give its description,
we write
$\mathfrak{g} = \mathfrak{u}_- + \mathfrak{l} + \mathfrak{u}_+$
for the Gelfand--Naimark decomposition,
and take a Cartan subalgebra $\mathfrak{j}$ of $\mathfrak{l}$ such that
$\mathfrak{l}^\tau$
contains
$\mathfrak{j}^\tau$
as a maximal abelian subspace (see \eqref{eqn:Vtau} for notation).
Let 
 $\Delta({\mathfrak{u}}_-^{-\tau}, {\mathfrak{j}}^{\tau})$ be the set of weights 
 of ${\mathfrak{u}}_-^{-\tau}$
 with respect to ${\mathfrak{j}}^{\tau}$.  
The roots $\alpha$ and $\beta$ are said to be 
 {\it{strongly orthogonal}}
 if neither $\alpha + \beta$ nor $\alpha - \beta$ is a root.  
We take a maximal set of strongly orthogonal roots 
$\{\nu_1, \cdots, \nu_k\}$ 
 in  $\Delta({\mathfrak{u}}_-^{-\tau}, {\mathfrak{j}}^{\tau})$ 
 inductively as follows:
 $\nu_{j}$ is the highest root among the elements
 in $\Delta({\mathfrak{u}}_-^{-\tau}, {\mathfrak{j}}^{\tau})$
 that are strongly orthogonal to $\nu_1, \cdots, \nu_{j-1}$
 ($1 \le j \le k-1$).   
The cardinality $k$ coincides with the split rank of the semisimple symmetric space $G_{\mathbb{R}}/G_{\mathbb{R}}^{\tau}$.

Then we have
\begin{theorem}
\label{thm:mfbr1}
Suppose that ${\mathfrak {p}}$ and $\tau$ are as in Theorem \ref{thm:PAbel}.  
Then, 
 for any sufficiently negative $\lambda$, 
 the generalized Verma module $\ensub {\mathfrak {g}}{\lambda}{\mathfrak{p}}$
 decomposes into a multiplicity-free direct sum of generalized Verma modules
 of ${\mathfrak {g}}^{\tau}$: 
\begin{equation}\label{eqn:mfbr}
  \ensub{\mathfrak {g}}{\lambda}{\mathfrak {p}}|_{\mathfrak {g}^{\tau}}
  \simeq
  \bigoplus_{\substack{a_1 \ge \cdots \ge a_l \ge 0 \\ a_1, \cdots, a_l \in {\mathbb{N}}}}
  \ensub{\mathfrak {g}^{\tau}}{\lambda|_{\mathfrak {j}^{\tau}}+\sum_{j=1}^l a_j \nu_j}{\mathfrak {p}^{\tau}}.  
\end{equation}
\end{theorem}

\begin{proof}[Proof of Theorem \ref{thm:mfbr1}]
Suppose that ${\mathfrak{p}}$ is a parabolic subalgebra
 such that its nilradical ${\mathfrak{u}}_+$ is abelian.  
Then $\mathfrak{p}$ is automatically a maximal parabolic subalgebra.
Further, it follows from  \cite{RRS92} that
 there exists a real form ${\mathfrak{g}}_{\mathbb{R}}$ of ${\mathfrak{g}}$
 such that $G_{\mathbb{R}}/(G_{\mathbb{R}} \cap P)$ is a Hermitian symmetric space
 of non-compact type,
 where $G_{\mathbb{R}}$ is the connected real form of $G={\operatorname{Int}}({\mathfrak{g}})$
 with Lie algebra ${\mathfrak{g}}_{\mathbb{R}}$.
The group 
$K_{\mathbb{R}} := G_{\mathbb{R}} \cap P$
is a maximal compact subgroup of $G_{\mathbb{R}}$,
and the complexification of its Lie algebra gives a Levi part,
denoted by $\mathfrak{l}$,
of $\mathfrak{p}$. 

Let $\theta$ be the involution of $\mathfrak{g}$ defined by
\[
\theta|_{\mathfrak{l}} = \operatorname{id},
\quad
\theta|_{\mathfrak{u}_- + \mathfrak{u}_+} = -\operatorname{id}.
\]
Then,
$\theta$ stabilizes $\mathfrak{g}_{\mathbb{R}}$ and $\mathfrak{p}$,
and the restriction
$\theta|_{\mathfrak{g}_{\mathbb{R}}}$
is a Cartan involution of the real semisimple Lie algebra
$\mathfrak{g}_{\mathbb{R}}$.
Since $\theta$ commutes with $\tau$,
$\tau\theta$ defines another involution of $\mathfrak{g}$.
We use the same symbol to denote its lift to the group $G$.
Then
$K_{\mathbb{R}}^{\tau\theta}
 = G_{\mathbb{R}}^{\tau\theta} \cap P$
is a maximal compact subgroup of 
$G_{\mathbb{R}}^{\tau\theta}$,
and has a complexified Lie algebra $\mathfrak{l}^\tau$.
Further,
$G_{\mathbb{R}}^{\tau\theta} / (G_{\mathbb{R}}^{\tau\theta} \cap P)
 = G_{\mathbb{R}}^{\tau\theta} / K_{\mathbb{R}}^{\tau\theta}$
becomes also a Hermitian symmetric space whose holomorphic tangent
space at the origin is identified with $\mathfrak{u}_-^{-\tau}$.
It then follows from W. Schmid \cite{xsch}
 that the symmetric algebra
$S(\mathfrak{u}_-^{-\tau})$
decomposes into the multiplicity-free sum of simple
$\mathfrak{l}^\tau$-modules as
\[
S(\mathfrak{u}_-^{-\tau}) \simeq \bigoplus_{\delta\in D} F'_\delta,
\]
where $\delta$ is the highest weight of
$F'_\delta$ and
\[
D := \{ \sum_{j=1}^k a_j \nu_j :
        a_1 \ge \dots \ge a_k \ge 0, \ 
        a_1, \dots, a_k \in \mathbb{N} \}.
\]

Applying Corollary \ref{cor:br},
we see that the identity \eqref{eqn:mfbr} holds in the Grothendieck group of
$\mathfrak{g}^\tau$-modules.
Finally,
let us show that the restriction
$\ensub{\mathfrak {g}}{\lambda}{\mathfrak {p}}|_{\mathfrak {g}^{\tau}}
$
decomposes as a direct sum of $\mathfrak{g}^\tau$-modules as given in
\eqref{eqn:mfbr} if $\lambda$ is
sufficiently negative.

For this,
let $\widetilde{G_{\mathbb{R}}}$ be the universal covering group of
$G_{\mathbb{R}}$,
and $\widetilde{K_{\mathbb{R}}}$ that of $K_{\mathbb{R}}$.
Then the generalized Verma module
$\ensub{\mathfrak{g}}{\lambda}{\mathfrak{p}}$
is isomorphic to the underlying
$(\mathfrak{g}, \widetilde{K_{\mathbb{R}}})$-module of a highest
weight representation of $\widetilde{G_{\mathbb{R}}}$
which is unitarizable if
$\langle \lambda,\alpha \rangle \ll 0$
for any 
$\alpha \in \Delta(\mathfrak{u}_+)$.
Hence, the identity \eqref{eqn:mfbr} in the Grothendieck group holds
as $\mathfrak{g}^\tau$-modules.
\end{proof}

\begin{remark}\label{rem:mfbr}
As we have seen in the above proof,
Theorems \ref{thm:PAbel} and \ref{thm:mfbr1} are equivalent to the
theorems on branching laws of unitary highest weight representations
of a real semisimple Lie group
$\widetilde{G_{\mathbb{R}}}$.
In the latter formulation,
the corresponding results were previously proved in \cite[Theorem B]{Ko07} by a
 geometric method based on reproducing kernels and `visible
actions' on complex manifolds \cite{Ko07b}.
See also \cite[Theorem 8.3]{Ko07} for explicit formulas.
\end{remark}

\subsection{Multiplicity-free pairs}
\label{subsec:MFn1}
Next, we consider multiplicity-free branching laws of the restriction 
$\ensub{\mathfrak{g}}{\lambda}{\mathfrak{p}}$
in the case where
$\mathfrak{p} = \mathfrak{b}$
(Borel subalgebra).
In general,
the `smaller' the parabolic subalgebra ${\mathfrak{p}}$ is, 
 the `larger' the generalized Verma module $\ensub{\mathfrak{g}}{\lambda}{\mathfrak{p}}$ becomes.  
Hence,
 we expect that the multiplicity-free property
 of the restriction $\ensub{\mathfrak{g}}{\lambda}{\mathfrak{b}}|_{\mathfrak{g}^{\tau}}$
 in the extreme case ${\mathfrak{p}}={\mathfrak{b}}$ 
should give the strongest constraints on the pair $({\mathfrak{g}}, {\mathfrak{g}}^{\tau})$.  
In this subsection,
 we determine for which symmetric pair $({\mathfrak{g}}, {\mathfrak{g}}^{\tau})$ 
 the restriction $\ensub{\mathfrak{g}}{\lambda}{\mathfrak{b}}|_{\mathfrak{g}^\tau}$
 is still multiplicity-free.  

Before stating a theorem,
we recall from Corollary \ref{cor:bdeco} that any simple
$\mathfrak{g}$-module in $\mathcal{O}$ contains at least one simple
$\mathfrak{g}^\tau$-module if and only if $G^\tau B$ is closed in $G$,
or equivalently, $\mathfrak{b}$ is $\tau$-stable.

\begin{theorem}
\label{thm:mfall}
Let ${\mathfrak{g}}$ be a complex simple Lie algebra,
 and $({\mathfrak{g}}, {\mathfrak{g}}^{\tau})$ a complex symmetric pair.  
Then the following three conditions are equivalent:
\begin{enumerate}
\item[{\rm{(i)}}]
$({\mathfrak{g}}, {\mathfrak{g}}^{\tau})$ is isomorphic to 
$
    ({\mathfrak{sl}}_{n+1}({\mathbb{C}}), {\mathfrak{gl}}_n({\mathbb{C}}))
$
 or
$
    ({\mathfrak{so}}_{n+1}({\mathbb{C}}), {\mathfrak{so}}_n({\mathbb{C}})).  
$

\item[{\rm{(ii)}}]
For any $\tau$-stable Borel subalgebra ${\mathfrak{b}}$,
 the restriction $\ensub{{\mathfrak{g}}}{\lambda}{{\mathfrak{b}}}|_{\mathfrak {g}^{\tau}}$
 is multiplicity-free as
$\mathfrak{g}^\tau$-modules for any generic $\lambda$.  

\item[{\rm{(iii)}}]
The restriction
$\ensub{\mathfrak{g}}{\lambda}{\mathfrak{b}}|_{\mathfrak{g}^\tau}$
is multiplicity-free as
$\mathfrak{g}^\tau$-modules for some $\lambda$ and some $\tau$-stable Borel
subalgebra $\mathfrak{b}$.
\end{enumerate}
\end{theorem}

\begin{proof}
[Proof of Theorem \ref{thm:mfall}]

(i) $\Longrightarrow$ (ii).\enspace
We shall give an explicit branching law of the restriction
$\ensub{\mathfrak{g}}{\lambda}{\mathfrak{b}}$
with respect to the symmetric
 $({\mathfrak{g}}, {\mathfrak{g}}^{\tau})$
 which is isomorphic to $({\mathfrak{sl}}_{n+1}({\mathbb{C}}), {\mathfrak{gl}}_{n}({\mathbb{C}}))$
 or $({\mathfrak{so}}_{n+1}({\mathbb{C}}), {\mathfrak{so}}_{n}({\mathbb{C}}))$
 in Subsections \ref{subsec:AA}--\ref{subsec:DB}.  

(ii) $\Longrightarrow$ (iii). \enspace
Obvious.

(iii) $\Longrightarrow$ (i). \enspace
We take a $\tau$-stable Levi decomposition 
 ${\mathfrak{b}}={\mathfrak{j}}+{\mathfrak{n}}$.  
Then, 
 it follows from Theorem \ref{thm:br}
 that $\ensub{\mathfrak {g}}{\lambda}{\mathfrak{b}}|_{\mathfrak {g}^{\tau}}$
 is multiplicity-free 
 only  if $S({\mathfrak {n}}^{-\tau})$ is multiplicity-free
 as a ${\mathfrak {j}}^{\tau}$-module.  
In turn,
 this happens only if
 the weights of ${\mathfrak {n}}^{-\tau}$
 are linearly independent over ${\mathbb{Q}}$, 
 which leads us to the following inequality
\[
     \dim {\mathfrak {n}}^{-\tau} \le \dim {\mathfrak {j}}^{\tau}, 
\]
or equivalently,
\begin{equation}
\label{eqn:dimrank}
\dim {\mathfrak {g}} -\dim {\mathfrak {g}}^{\tau}
  \le {\operatorname{{rank}}}\,{\mathfrak {g}}+ \operatorname{rank}{\mathfrak {g}}^{\tau}.  
\end{equation}
In view of the classification of complex symmetric pairs
 $({\mathfrak{g}}, {\mathfrak{g}}^{\tau})$ 
with $\mathfrak{g}$ simple, 
 the inequality \eqref{eqn:dimrank} holds 
 only if $({\mathfrak {g}}, {\mathfrak {g}}^{\tau})$
 is isomorphic to $({\mathfrak {sl}}_{n+1}({\mathbb{C}}), {\mathfrak {gl}}_n({\mathbb{C}}))$
 or $({\mathfrak {so}}_{n+1}({\mathbb{C}}), {\mathfrak {so}}_{n}({\mathbb{C}}))$.  
\end{proof}

In Subsections \ref{subsec:AA}--\ref{subsec:DB},
we shall fix a Borel subalgebra $\mathfrak{b}$ of $\mathfrak{g}$
and consider $B$-conjugacy classes of
involutions $\tau$
instead of considering $G^\tau$-conjugacy classes of Borel subalgebras
by fixing $\tau$.
With this convention,
we shall use the abbreviation
$\en{\mathfrak{g}}{\lambda}$ for
$\ensub{\mathfrak{g}}{\lambda}{\mathfrak{b}}$.

\subsection{Branching laws for 
$\protect\mathfrak{gl}_{n+1} \downarrow \mathfrak{gl}_n$}
\label{subsec:AA}

Let ${\mathfrak{g}}:={\mathfrak{gl}}_{n+1}({\mathbb{C}})$
 and ${\mathfrak{g}}':={\mathfrak{gl}}_1({\mathbb{C}}) \oplus {\mathfrak{gl}}_n({\mathbb{C}})$.  
We observe that there are $(n+1)$ closed $GL_n({\mathbb{C}})$-orbits on the full flag variety of $GL_{n+1}({\mathbb{C}})$.  
Correspondingly, 
 there are essentially $n+1$ different settings for discretely
decomposable restrictions of
the Verma module $\en{\mathfrak{g}}{\lambda}$ to $\mathfrak{g}'$
 by Theorems \ref{thm:orbit} and \ref{thm:decoc}.

In order to fix notation, 
 let ${\mathfrak {b}}=
{\mathfrak{j}}+{\mathfrak{n}}_+$ be the standard Borel subalgebra
 of consisting of upper triangular matrices in
${\mathfrak {g}}$,
and $\mathfrak{j}$ the Cartan subalgebra consisting of diagonal matrices.
For $1 \le l \le n+1$,
 we realize ${\mathfrak {g}}'$
 as a subalgebra of ${\mathfrak {g}}$ 
 by letting $\iota_l({\mathfrak {g}}')$ be the centralizer
 of the matrix unit $E_{ll}$.   
For $k=(k_1, \cdots, \widehat{k_l}, \cdots, k_{n+1}) \in {\mathbb{N}}^n$, 
 we set 
\[
  \operatorname{ind}_l k:= k_1 + \cdots + k_{l-1} - k_{l+1} -\cdots - k_{n+1}.  
\]

In what follows,
 $\boxtimes$ denotes the outer tensor product representation of the direct product
of Lie algebras.  

\begin{theorem}
[$A_n \downarrow A_{n-1}$]
\label{thm:AA}
Suppose $\lambda_i - \lambda_j \notin \mathbb{Z}$
for any distinct $i, j$ in $\{1, \cdots, \widehat l, \cdots, n+1\}$.  
Then the restriction of the Verma module of\/
$\mathfrak{g}$
decomposes into a multiplicity-free direct sum of
simple Verma modules of\/
$\mathfrak{g}'$.
\begin{multline}
\label{eqn:AA}
  \en{{\mathfrak{gl}}_{n+1}}{\lambda}|_{\iota_l ({\mathfrak{gl}_1}\oplus {\mathfrak{gl}}_n)}
\\
  \simeq 
   \bigoplus_{k \in {\mathbb{N}}^n}
  {\mathbb{C}}_{\lambda_l + \operatorname{ind}_l k}
   \boxtimes
   \en{\mathfrak {gl}_n}
   {\lambda_1 - k_1, \cdots, \lambda_{l-1}- k_{l-1}, 
   \lambda_{l+1} + k_{l+1}, \cdots, \lambda_{n+1} + k_{n+1}}.  
\end{multline}
\end{theorem}

\medskip

\begin{proof}
We fix
$l$ ($1 \le l \le n+1$)
once and for all.
Let $\tau \equiv \tau_l$ be the involution of $\mathfrak{g}$ such that
$\mathfrak{g}^\tau = \iota_l(\mathfrak{g}')$.
With our choice of $\mathfrak{j}$,
we have
$\mathfrak{j}^\tau = \mathfrak{j} \simeq \mathbb{C}^{n+1}$,
and the set of characters of $\mathfrak{j}^\tau$ is identified with
$\mathbb{C}^{n+1}$.
We apply Corollary \ref{cor:br}
 to the $\mathfrak{j}^{\tau}$-module $\mathfrak{n}_-^{-\tau}$: 
\[
     {\mathfrak {n}}_-^{-\tau} 
    = \bigoplus_{i=1}^{l-1} \mathfrak{g}_{-e_i + e_l}
       \oplus
       \bigoplus_{i=l+1}^{n+1} \mathfrak{g}_{-e_l + e_j}.  
\]
Extending this to the symmetric algebra $S({\mathfrak {n}}_-^{-\tau})$,
 we have ${\mathfrak{j}}^{\tau}$-isomorphism: 
\[
     S({\mathfrak {n}}_-^{-\tau}) \simeq
     \bigoplus_{k \in {\mathbb{N}}^n}
    (-k_1, \dots, -k_{l-1}, \operatorname{ind}_l k, k_{l+1}, \cdots, k_{n+1}).  
\]
Therefore, 
 the identity \eqref{eqn:AA} holds
 in the Grothendieck group
 by Corollary \ref{cor:br}.  

Since $\lambda_i - \lambda_j \not \in {\mathbb{Z}}$
 for any $i, j$, 
 the Verma modules appearing in the right-hand side of \eqref{eqn:AA}
 have distinct infinitesimal characters.  
Therefore, 
 there is no extension among these representations.  
Hence \eqref{eqn:AA} is a direct sum decomposition.
\end{proof}

\subsection{Branching laws for
$\protect\mathfrak{so}(2n+1) \downarrow \mathfrak{so}(2n)$}
\label{subsec:BD}
Let $\mathfrak{g} = \mathfrak{so}_{2n+1}(\mathbb{C})$, 
 $\mathfrak{g}' = \mathfrak{so}_{2n}(\mathbb{C})$
 and $G'$ be the connected subgroup of $G=\operatorname{Int}({\mathfrak{g}})$
 with Lie algebra ${\mathfrak{g}}'$.  
Then there are two closed $G'$-orbits
 on the full flag variety $G/B$, 
 which are conjugate to each other
 by an element of the normalizer $N_G({\mathfrak{g}}')$.  
Thus  
 it follows from Theorem \ref{thm:orbit}
 that there is essentially the unique triple 
$
 ({\mathfrak{g}}, {\mathfrak{g}}',{\mathfrak{b}})
$
 satisfying the equivalent conditions
 of Theorem \ref{thm:decoc}.

To fix notation, 
 we may and do assume  
 that ${\mathfrak{g}}' \cap {\mathfrak{b}}$
 contains a Cartan subalgebra ${\mathfrak{j}}$ of ${\mathfrak{g}}$
 and that  
\begin{align*}
\Delta^+(\mathfrak{g},\mathfrak{j})
=& \{ e_i \pm e_j : 1 \le i < j \le n \} \cup \{ e_i : 1 \le i \le n \}, 
\\
\Delta^+(\mathfrak{g}',\mathfrak{j})
=& \{ e_i \pm e_j : 1 \le i < j \le n \}.   
\end{align*}
\begin{theorem}
[$B_n \downarrow D_n$]
\label{thm:BD}
Suppose $\lambda_i \pm \lambda_j \notin \mathbb{Z}$
for any $1 \le i < j \le n$.
\begin{equation}
\label{eqn:BD}
\en{\mathfrak{so}_{2n+1}}{\lambda} |_{\mathfrak{so}_{2n}}
= \bigoplus_{k\in\mathbb{N}^n} \en{\mathfrak{so}_{2n}}{\lambda- k}.
\end{equation}
\end{theorem}

\begin{proof}
[Proof of Theorem \ref{thm:BD}]
Let $\tau$ be the involution of ${\mathfrak{g}}$
 such that ${\mathfrak{g}}'={\mathfrak{g}}^{\tau}$.  
Applying Corollary \ref{cor:br}
 to the ${\mathfrak{j}}$-module:
\begin{equation*}
  S(\mathfrak{n}_-^{-\tau}) =  S(\bigoplus_{i=1}^{n} \mathfrak{g}_{-e_i})
  \simeq \bigoplus_{k \in {\mathbb{N}}^n} (-k_1, \cdots, -k_n), 
\end{equation*}  
 we get \eqref{eqn:BD} 
 in the Grothendieck group.  
The assumption $\lambda_i \pm \lambda_j \not \in {\mathbb{Z}}$
 assures that every summand in \eqref{eqn:BD} is simple.  
Further,
 there is no extension among $\en{\mathfrak{so}_{2n}}{\lambda -k}$
 because they have a distinct ${\mathfrak{Z}}({\mathfrak{g}'})$-infinitesimal characters.   
\end{proof}

\subsection{Branching laws for
$\mathfrak{so}(2n+2) \downarrow \mathfrak{so}(2n+1)$}
\label{subsec:DB}
Let ${\mathfrak{g}}={\mathfrak{so}}_{2n+2}(\mathbb{C})$
 and ${\mathfrak{g}}'={\mathfrak{so}}_{2n+1}(\mathbb{C})$.  
Then there exists a unique closed $G'$-orbit 
 on the full flag variety $G/B$.  
To fix notation,
 we suppose that our Borel subalgebra
$
    {\mathfrak{b}}={\mathfrak{j}}+{\mathfrak{n}}
$
 is defined by the positive system
$$
\Delta^+(\mathfrak{g},\mathfrak{j})
= \{ e_i \pm e_j : 1 \le i < j \le n+1 \},  
$$
 and that 
$
 {\mathfrak{j}}':={\mathfrak{j}} \cap {\mathfrak{g}}'
$
 is given by 
$
\{
     H \in  {\mathfrak{j}}: e_{n+1}(H) =0
\}.  
$
Then ${\mathfrak{b}}':={\mathfrak{b}} \cap {\mathfrak{g}}'$
 is a Borel subalgebra of ${\mathfrak{g}}'$
 given by a positive system 
\[
  \Delta^+({\mathfrak{g}}', {\mathfrak{j}}')
  =\{ e_i \pm e_j: 1 \le i < j \le n\} \cup \{ e_i : 1 \le i \le n\}.  
\]
\begin{theorem}
[$D_{n+1} \downarrow B_n$]
\label{thm:DB}
Suppose $\lambda_i \pm \lambda_j \not\in {\mathbb{Z}}$
 for any $1 \le i < j \le n$.  
We set $\lambda:=(\lambda_1, \cdots, \lambda_n)$.  
Then
\begin{equation}
\label{eqn:DB}
\en{{\mathfrak{so}}_{2n+2}}{\lambda,\lambda_{n+1}}|_{{\mathfrak{so}}_{2n+1}}
  \simeq \bigoplus_{k \in {\mathbb{N}}^n}
  \en{{\mathfrak{so}}_{2n+1}}{\lambda -k }.  
\end{equation}
\end{theorem}
\begin{proof}
Let $\tau$ be the defining involution of ${\mathfrak{g}}'={\mathfrak{so}}_{2n+1}({\mathbb{C}})$.  
Then 
\begin{align*}
  {\mathfrak{n}}_-^{-\tau} 
  =& \bigoplus_{i=1}^n ({\mathfrak{g}}_{-e_i + e_{n+1}} + {\mathfrak{g}}_{-e_i -e_{n+1}})^{-\tau}, 
\\
\intertext{and hence we have an isomorphism }
  S({\mathfrak{n}}_-^{-\tau}) \simeq& \bigoplus_{k \in {\mathbb{N}^n}} (-k_1, \cdots, -k_n)
\end{align*}
as ${\mathfrak{j}}'$-modules.  
Therefore, 
 \eqref{eqn:DB} follows from Corollary \ref{cor:br}.  
\end{proof}

\bigskip

\noindent
\textsc{\small
IMPU, 
Graduate School of Mathematical Sciences,
The University of Tokyo,
3-8-1 Komaba, Tokyo, 153-8914 Japan.
%;
%Institut des Hautes \'{E}tudes Scientifiques,
%Bures-sur-Yvette, France (current address).
}

\end{document}